\documentclass{amsart}

\usepackage{amsfonts,amssymb,verbatim,amsmath,amsthm,latexsym,textcomp,amscd}
\usepackage{latexsym,amsfonts,amssymb,epsfig,verbatim}
\usepackage{amsmath,amsthm,amssymb,latexsym,graphics,textcomp}
\usepackage{graphicx}
\usepackage{color}
\usepackage{url}

\begin{document}

\newtheorem{theorem}{Theorem}[section]
\newtheorem{prop}[theorem]{Proposition}
\newtheorem{lemma}[theorem]{Lemma}
\newtheorem{cor}[theorem]{Corollary}
\newtheorem{definition}[theorem]{Definition}
\newtheorem{conj}[theorem]{Conjecture}
\newtheorem{claim}[theorem]{Claim}
\newtheorem{defth}[theorem]{Definition-Theorem}
\newtheorem{obs}[theorem]{Observation}
\newtheorem{rmk}[theorem]{Remark}
\newtheorem{qn}[theorem]{Question}

\newcommand{\hhat}{\widehat}
\newcommand{\boundary}{\partial}
\newcommand{\C}{{\mathbb C}}
\newcommand{\integers}{{\mathbb Z}}
\newcommand{\natls}{{\mathbb N}}
\newcommand{\ratls}{{\mathbb Q}}
\newcommand{\reals}{{\mathbb R}}
\newcommand{\proj}{{\mathbb P}}
\newcommand{\lhp}{{\mathbb L}}
\newcommand{\tube}{{\mathbb T}}
\newcommand{\cusp}{{\mathbb P}}
\newcommand\AAA{{\mathcal A}}
\newcommand\BB{{\mathcal B}}
\newcommand\CC{{\mathcal C}}
\newcommand\DD{{\mathcal D}}
\newcommand\EE{{\mathcal E}}
\newcommand\FF{{\mathcal F}}
\newcommand\GG{{\mathcal G}}
\newcommand\HH{{\mathcal H}}
\newcommand\II{{\mathcal I}}
\newcommand\JJ{{\mathcal J}}
\newcommand\KK{{\mathcal K}}
\newcommand\LL{{\mathcal L}}
\newcommand\MM{{\mathcal M}}
\newcommand\NN{{\mathcal N}}
\newcommand\OO{{\mathcal O}}
\newcommand\PP{{\mathcal P}}
\newcommand\QQ{{\mathcal Q}}
\newcommand\RR{{\mathcal R}}
\newcommand\SSS{{\mathcal S}}
\newcommand\TT{{\mathcal T}}
\newcommand\UU{{\mathcal U}}
\newcommand\VV{{\mathcal V}}
\newcommand\WW{{\mathcal W}}
\newcommand\XX{{\mathcal X}}
\newcommand\YY{{\mathcal Y}}
\newcommand\ZZ{{\mathcal Z}}
\newcommand\CH{{\CC\Hyp}}
\newcommand\MF{{\MM\FF}}
\newcommand\PMF{{\PP\kern-2pt\MM\FF}}
\newcommand\ML{{\MM\LL}}
\newcommand\PML{{\PP\kern-2pt\MM\LL}}
\newcommand\GL{{\GG\LL}}
\newcommand\Pol{{\mathcal P}}
\newcommand\half{{\textstyle{\frac12}}}
\newcommand\Half{{\frac12}}
\newcommand\Mod{\operatorname{Mod}}
\newcommand\Area{\operatorname{Area}}
\newcommand\ep{\epsilon}
\newcommand\Hypat{\widehat}
\newcommand\Proj{{\mathbf P}}
\newcommand\U{{\mathbf U}}
 \newcommand\Hyp{{\mathbf H}}
\newcommand\D{{\mathbf D}}
\newcommand\Z{{\mathbb Z}}
\newcommand\R{{\mathbb R}}
\newcommand\Q{{\mathbb Q}}
\newcommand\E{{\mathbb E}}
\newcommand\til{\widetilde}
\newcommand\length{\operatorname{length}}
\newcommand\tr{\operatorname{tr}}
\newcommand\gesim{\succ}
\newcommand\lesim{\prec}
\newcommand\simle{\lesim}
\newcommand\simge{\gesim}
\newcommand{\simmult}{\asymp}
\newcommand{\simadd}{\mathrel{\overset{\text{\tiny $+$}}{\sim}}}
\newcommand{\ssm}{\setminus}
\newcommand{\diam}{\operatorname{diam}}
\newcommand{\pair}[1]{\langle #1\rangle}
\newcommand{\T}{{\mathbf T}}
\newcommand{\inj}{\operatorname{inj}}
\newcommand{\pleat}{\operatorname{\mathbf{pleat}}}
\newcommand{\short}{\operatorname{\mathbf{short}}}
\newcommand{\vertices}{\operatorname{vert}}
\newcommand{\collar}{\operatorname{\mathbf{collar}}}
\newcommand{\bcollar}{\operatorname{\overline{\mathbf{collar}}}}
\newcommand{\I}{{\mathbf I}}
\newcommand{\tprec}{\prec_t}
\newcommand{\fprec}{\prec_f}
\newcommand{\bprec}{\prec_b}
\newcommand{\pprec}{\prec_p}
\newcommand{\ppreceq}{\preceq_p}
\newcommand{\sprec}{\prec_s}
\newcommand{\cpreceq}{\preceq_c}
\newcommand{\cprec}{\prec_c}
\newcommand{\topprec}{\prec_{\rm top}}
\newcommand{\Topprec}{\prec_{\rm TOP}}
\newcommand{\fsub}{\mathrel{\scriptstyle\searrow}}
\newcommand{\bsub}{\mathrel{\scriptstyle\swarrow}}
\newcommand{\fsubd}{\mathrel{{\scriptstyle\searrow}\kern-1ex^d\kern0.5ex}}
\newcommand{\bsubd}{\mathrel{{\scriptstyle\swarrow}\kern-1.6ex^d\kern0.8ex}}
\newcommand{\fsubeq}{\mathrel{\raise-.7ex\hbox{$\overset{\searrow}{=}$}}}
\newcommand{\bsubeq}{\mathrel{\raise-.7ex\hbox{$\overset{\swarrow}{=}$}}}
\newcommand{\tw}{\operatorname{tw}}
\newcommand{\base}{\operatorname{base}}
\newcommand{\trans}{\operatorname{trans}}
\newcommand{\rest}{|_}
\newcommand{\bbar}{\overline}
\newcommand{\UML}{\operatorname{\UU\MM\LL}}
\newcommand{\EL}{\mathcal{EL}}
\newcommand{\tsum}{\sideset{}{'}\sum}
\newcommand{\tsh}[1]{\left\{\kern-.9ex\left\{#1\right\}\kern-.9ex\right\}}
\newcommand{\Tsh}[2]{\tsh{#2}_{#1}}
\newcommand{\qeq}{\mathrel{\approx}}
\newcommand{\Qeq}[1]{\mathrel{\approx_{#1}}}
\newcommand{\qle}{\lesssim}
\newcommand{\Qle}[1]{\mathrel{\lesssim_{#1}}}
\newcommand{\simp}{\operatorname{simp}}
\newcommand{\vsucc}{\operatorname{succ}}
\newcommand{\vpred}{\operatorname{pred}}
\newcommand\fhalf[1]{\overrightarrow {#1}}
\newcommand\bhalf[1]{\overleftarrow {#1}}
\newcommand\sleft{_{\text{left}}}
\newcommand\sright{_{\text{right}}}
\newcommand\sbtop{_{\text{top}}}
\newcommand\sbot{_{\text{bot}}}
\newcommand\sll{_{\mathbf l}}
\newcommand\srr{_{\mathbf r}}
\newcommand\geod{\operatorname{\mathbf g}}
\newcommand\mtorus[1]{\boundary U(#1)}
\newcommand\A{\mathbf A}
\newcommand\Aleft[1]{\A\sleft(#1)}
\newcommand\Aright[1]{\A\sright(#1)}
\newcommand\Atop[1]{\A\sbtop(#1)}
\newcommand\Abot[1]{\A\sbot(#1)}
\newcommand\boundvert{{\boundary_{||}}}
\newcommand\storus[1]{U(#1)}
\newcommand\Momega{\omega_M}
\newcommand\nomega{\omega_\nu}
\newcommand\twist{\operatorname{tw}}
\newcommand\modl{M_\nu}
\newcommand\MT{{\mathbb T}}
\newcommand\Teich{{\mathcal T}}
\renewcommand{\Re}{\operatorname{Re}}
\renewcommand{\Im}{\operatorname{Im}}

\title{Cannon-Thurston Maps,
 i-bounded Geometry and a Theorem of McMullen}

\author{Mahan Mj}
\address{RKM Vivekananda University}
 \date{} 

\thanks{Research partly supported by UGC Major Research Project}

\begin{abstract}
The notion of  {\it i-bounded geometry}  generalises
 simultaneously 
{\it bounded geometry} and the geometry of punctured torus Kleinian
 groups.  We
show that the 
limit set of a surface Kleinian group of i-bounded geometry
 is locally connected by constructing a natural
Cannon-Thurston map.

\smallskip

\begin{center}


\end{center}

\end{abstract}

\maketitle

\tableofcontents

\section{Introduction} In \cite{mahan-split} we prove the existence of Cannon-Thurston maps 
for arbitrary surface Kleinian groups without accidental parabolics.  The proof proceeds by constructing a coarse
model geometry, called {\it split geometry},
satisfied by all associated hyperbolic 3-manifolds. Our starting point in \cite{mahan-split} is a model
geometry constructed by Minsky in \cite{minsky-elc1} and we proceed by forgetting some of the
finer structure in \cite{minsky-elc1} to establish that all  surface Kleinian groups have 
associated hyperbolic 3-manifolds of split geometry. In \cite{mahan-red}, \cite{mahan-elct}, \cite{mahan-elct2} and \cite{mahan-kl} we 
completed the programme of proving the existence of Cannon-Thurston maps 
for arbitrary finitely generated Kleinian groups and describing point pre-images in terms of ending laminations.

The purpose of the present paper is to give an exposition 
of the existence of Cannon-Thurston maps for surface Kleinian groups without accidental parabolics satisfying
a more restrictive model geometry called {\bf i-bounded geometry} satisfied for instance by all punctured torus
Kleinian groups. This gives a new proof of a result of McMullen \cite{ctm-locconn}.

The main pre-requisites for understanding the present paper are: \\
1) Generalities on hyperbolic metric
spaces  in the sense
of Gromov \cite{gromov-hypgps} \cite{CDP}, \cite{GhH}, especially boundary theory in terms of asymptote classes of geodesics. \\
2) The theory of simply and doubly degenerate Kleinian groups (Chapter 8 of \cite{thurstonnotes}) \\
3) Relative hyperbolicity and electric geometry \cite{farb-relhyp}, \cite{gromov-hypgps}, \cite{bowditch-relhyp}.

\smallskip

A similar exposition in the special case of bounded geometry surface Kleinian groups had been given by the author in \cite{brahma-bddgeo}.
In \cite{brahma-amalgeo} we give an exposition of more general model geometries leading up to split geometry used in \cite{mahan-split}.

\subsection{Statement of Results}

The main theorem of this paper is:

\smallskip

{\bf Theorem \ref{crucial-punct}}: {\it
Let $M^h$ be a hyperbolic 
3 manifold of {\bf i-bounded geometry} homeomorphic to $S^h \times J$ (for $J = [0,
  \infty ) $ or $( - \infty , \infty )$), where $S^h$ is a hyperbolic surface of finite area.
 Let $i: S^h \rightarrow M^h$ be a type-preserving (i.e. taking
 parabolics to parabolics) homotopy equivalence.  Then the inclusion
  $\tilde{i} : \widetilde{S^h} \rightarrow \widetilde{M^h}$ extends
 continuously 
  to a map 
  $\hat{i} : \widehat{S^h} \rightarrow \widehat{M^h}$. Hence the limit set
  of $\widetilde{S^h}$ is locally connected.}

\smallskip

The notion of {\it i-bounded geometry}  generalises
 simultaneously 
{\it bounded geometry} and the geometry of punctured torus Kleinian
 groups. In particular, since punctured torus groups have {\em
 i-bounded geometry} by a result of Minsky \cite{minsky-torus}, we
 have a new proof of the following Theorem of McMullen
 \cite{ctm-locconn} as a consequence: \\

\smallskip

\noindent {\bf Theorem :} (McMullen \cite{ctm-locconn} ) {\it
Let $M^h$ be a hyperbolic 
3 manifold homeomorphic to $S^h \times J$ (for $J = [0,
  \infty ) $ or $( - \infty , \infty )$), where $S^h$ is a punctured
  torus.
 Let $i: S^h \rightarrow M^h$ be a type-preserving (i.e. taking
 parabolics to parabolics) homotopy equivalence.  Then the inclusion
  $\tilde{i} : \widetilde{S^h} \rightarrow \widetilde{M^h}$ extends
 continuously 
  to a map 
  $\hat{i} : \widehat{S^h} \rightarrow \widehat{M^h}$. Hence the limit set
  of $\widetilde{S^h}$ is locally connected.}

\smallskip

{\em i-bounded geometry}
 can roughly be described as bounded geometry away from Margulis
 tubes. But this description is a little ambiguous. More precisely, we
 start with a collection of (uniformly) bounded geometry blocks $S
 \times I$ glued end to end. Next, for some blocks a curve is selected
 such that its representative on the lower end of the block has
(uniformly)  bounded length. Hyperbolic Dehn surgery is then performed
 along the geodesic representative within the block. Precise
 definitions will be given in Section 2.2.

We describe below a collection of  examples of manifolds of {\it
  i-bounded geometry} for which Theorem \ref{crucial-punct} is known:

\noindent 1)
 The cover corresponding to the fiber subgroup of a closed
  hyperbolic 3-manifold fibering over the circle (Cannon and Thurston
  \cite{CT}). \\
\noindent 2) Hyperbolic 3 manifolds of bounded geometry, which correspond to
  simply or doubly degenerate Kleinian groups isomorphic to closed
  surface groups (Minsky \cite{minsky-jams}). (See also   Section 4.3 of
 \cite{mitra-trees}.)\\
\noindent 3) Hyperbolic 3 manifolds of bounded geometry, arising from
  simply or doubly degenerate Kleinian groups corresponding to punctured
  surface groups (Bowditch \cite{bowditch-ct}). (See also \cite{brahma-pared})\\
\noindent 4) Punctured torus Kleinian groups (McMullen
  \cite{ctm-locconn}). \\

\subsection{Cannon-Thurston Maps and i-bounded geometry}

Let $S$ be a hyperbolic surface of finite area and  let $\rho ( \pi_1
(S)) = H \subset PSl_2 ( {\Bbb{C}} ) $ = Isom ( ${\Bbb{H}}^3$) be a
representation, such that the quotient hyperbolic 3-manifold
$M = {{\Bbb{H}}^3}/H$  is simply degenerate.  
 Let $\widetilde S$ and $\widetilde M$ denote the universal
covers of $S$ and $M$ respectively. Then $\widetilde S$ and $\widetilde M$
can be identified with
 ${\Bbb{H}}^2$ and ${\Bbb{H}}^3$ respectively. There exists a natural
inclusion $i: \widetilde{S} \rightarrow \widetilde{M}$. Now let
${{\Bbb{D}}^2}={\Bbb{H}}^2\cup{\Bbb{S}}^1_\infty$ and 
${{\Bbb{D}}^3}={\Bbb{H}}^3\cup{\Bbb{S}}^2_\infty$
denote the standard compactifications. The local connectivity of the
limit set of $\widetilde S$ is equivalent to the existence of a
continuous extension (a {\bf Cannon-Thurston map}) $\hat{i} : 
{{\Bbb{D}}^2} \rightarrow {{\Bbb{D}}^3}$.

A word about the term {\bf i-bounded geometry}. 
In the construction of a general model
manifold (Section 9 of \cite{minsky-elc1}),
as a step towards the  resolution of the Ending
Lamination Conjecture, Minsky describes certain (complex) meridian
 coefficients which encode the complex structure for boundary torii of
Margulis tubes. The uniform boundedness of these coefficients
corresponds to {\it bounded geometry}. The
manifolds that we discuss in this paper correspond to those which have
a uniform bound on the {\bf imaginary}
 part of these coefficients. Hence the term
{\it i-bounded geometry}. Clearly, manifolds of {\it bounded geometry}
 have {\it i-bounded geometry}. In
\cite{minsky-torus}, Minsky further 
showed that punctured torus groups (and four-holed sphere groups) have
{\it i-bounded geometry}. Roughly speaking, the number of twists gives
the real part and the number of vertical annulii gives the imaginary
part of the coefficients. Hence, in a manifold of {\it i-bounded
  geometry}, an arbitrarily large number of twists are allowed for
each Margulis tube, but only a uniformly bounded number of vertical annulii.

As in \cite{mitra-ct}, \cite{mitra-trees} and \cite{brahma-pared}, our
proof proceeds  by constructing a {\it ladder-like} set $B_\lambda \subset
\widetilde{M}$  from a geodesic segment $\lambda \subset
\widetilde{S}$ and then a retraction $\Pi_\lambda$ of $\widetilde{M}$
onto $B_\lambda$. We modify this construction in this paper and
restrict our attention to one block, i.e. a copy of $\widetilde{S}
\times I$ minus certain  neighborhoods of geodesics and  cusps and
equip it with a model {\em pseudometric} which is zero along lifts of a
simple closed geodesic.

To prevent cluttering, we
restrict ourselves to closed surfaces first, and then indicate the
modifications necessary for punctured surfaces.

\section{Preliminaries}

\subsection{Hyperbolic Metric Spaces}

We start off with some preliminaries about hyperbolic metric
spaces  in the sense
of Gromov \cite{gromov-hypgps} \cite{CDP}, \cite{GhH}. Let $(X,d)$
be a hyperbolic metric space. The 
{\bf Gromov boundary} of 
 $X$, denoted by $\partial{X}$,
is the collection of asymptote classes of geodesic rays.

 A subset $Z$ of $X$ is said to be 
{\bf $k$-quasiconvex}
 if any geodesic joining points of  $ Z$ lies in a $k$-neighborhood of $Z$.
A subset $Z$ is {\bf quasiconvex} if it is $k$-quasiconvex for some
$k$. 

A map $f$ from one metric space $(Y,{d_Y})$ into another metric space 
$(Z,{d_Z})$ is said to be
 a {\bf $(K,\epsilon)$-quasi-isometric embedding} if
 
\begin{center}
${\frac{1}{K}}({d_Y}({y_1},{y_2}))-\epsilon\leq{d_Z}(f({y_1}),f({y_2}))\leq{K}{d_Y}({y_1},{y_2})+\epsilon$
\end{center}
If  $f$ is a quasi-isometric embedding, 
 and every point of $Z$ lies at a uniformly bounded distance
from some $f(y)$ then $f$ is said to be a {\bf quasi-isometry}.
A $(K,{\epsilon})$-quasi-isometric embedding that is a quasi-isometry
will be called a $(K,{\epsilon})$-quasi-isometry.

A {\bf $(K,\epsilon)$-quasigeodesic}
 is a $(K,\epsilon)$-quasi-isometric embedding
of
a closed interval in $\Bbb{R}$. A $(K,K)$-quasigeodesic will also be called
a $K$-quasigeodesic.

Let $(X,{d_X})$ be a hyperbolic metric space and $Y$ be a subspace that is
hyperbolic with the inherited path metric $d_Y$.
By 
adjoining the Gromov boundaries $\partial{X}$ and $\partial{Y}$
 to $X$ and $Y$, one obtains their compactifications
$\widehat{X}$ and $\widehat{Y}$ respectively.

Let $ i :Y \rightarrow X$ denote inclusion.

{\bf Definition:}   Let $X$ and $Y$ be hyperbolic metric spaces and
$i : Y \rightarrow X$ be an embedding. 
 A {\bf Cannon-Thurston map} $\hat{i}$  from $\widehat{Y}$ to
 $\widehat{X}$ is a continuous extension of $i$.

The following  lemma (Lemma 2.1 of \cite{mitra-ct})
 says that a Cannon-Thurston map exists
if for all $M > 0$ and $y \in Y$, there exists $N > 0$ such that if $\lambda$
lies outside an $N$ ball around $y$ in $Y$ then
any geodesic in $X$ joining the end-points of $\lambda$ lies
outside the $M$ ball around $i(y)$ in $X$.
For convenience of use later on, we state this somewhat
differently.

\begin{lemma}
A Cannon-Thurston map from $\widehat{Y}$ to $\widehat{X}$
 exists if  the following condition is satisfied:

Given ${y_0}\in{Y}$, there exists a non-negative function  $M(N)$, such that 
 $M(N)\rightarrow\infty$ as $N\rightarrow\infty$ and for all geodesic segments
 $\lambda$  lying outside an $N$-ball
around ${y_0}\in{Y}$  any geodesic segment in $\Gamma_G$ joining
the end-points of $i(\lambda)$ lies outside the $M(N)$-ball around 
$i({y_0})\in{X}$.

\label{contlemma}
\end{lemma}

\smallskip

The above result can be interpreted as saying that a Cannon-Thurston map 
exists if the space of geodesic segments in $Y$ embeds properly in the
space of geodesic segments in $X$.

\subsection{i-bounded Geometry}

We start with a hyperbolic surface $S^h$ with or without punctures. The
hyperbolic structure is arbitrary, but it is important that a choice
be made. $S$ will denote $S^h$ minus a small enough neighborhood of
the cusps.

  Fix a {\em finite}
collection $\mathcal{C}$ of (geodesic representatives
of) simple closed curves on
$S$. $N_\epsilon ( \sigma )$ will denote the $\epsilon$-neghborhood of
a geodesic $\sigma \in \mathcal{C}$. 

 $N_\epsilon ( \sigma_i )$ will denote an $\epsilon$ neighborhood of
$\sigma_i \subset S^h$ for some $\sigma_i \in \mathcal{C}$. $\epsilon$
and the neighborhood of the cusps 
in $S^h$ are chosen small enough so
that \\
$\bullet 1$  $N_\epsilon ( \sigma_i )$ is at least a distance of $\epsilon$
from the cusps. \\
$\bullet 2$ No two lifts of  $N_\epsilon ( \sigma_i )$ to the
universal cover $\widetilde{S^h}$ intersect. \\

Note that $S = S^h$ if $S$ has
no cusps. Restrict the metric on $S^h$ to $S$ and equip $S$ with the
resultant {\em path-metric}.

\smallskip

\noindent {\bf The Thin Building Block}\\
For the construction of a thin block, $I$ will denote the closed
interval $[0,3]$. 
Now put a product metric structure on $S \times I$, which restricts to
the path-metric on $S$ for each slice $S \times {a}, a \in I $ and 
the Euclidean metric on the $I$-factor. Let
$B_i^c$ denote $(S \times I - N_{\epsilon} ( \sigma_i ) \times
[1,2]$. 
Equip $B_i^c$ with the path-metric.

For each resultant torus component of
the boundary of $B_i^c$, perform Dehn filling on some $(1,n_i)$ curve,
which goes $n_i$ times around the meridian and once round the
longitude.  $n_i$ will be called the {\bf twist coefficient}.
The metric on the solid torus $\Theta_i$ glued in is arranged in such
a way that it is isometric to the quotient of a neighborhood of a
bi-infinite hyperbolic geodesic by a hyperbolic isometry. Further, the
$(1, n_i )$-curve is required to bound a totally geodesic hyperbolic
disk. In fact, we might as well foliate the boundary of $\Theta_i$ by
translates (under hyperbolic isometries) of the meridian, and demand
that each bounds a totally geodesic disk. Since there is
no canonical way to smooth out the resulting metric, we leave it as
such. $\Theta_i$ equipped with this metric will be called a {\bf
  Margulis tube} in keeping with the analogy from hyperbolic space.

The resulting copy of $S \times I$ obtained, equipped with the metric just
described, is called a {\bf thin building block} and is denoted by $B_i$. 

\smallskip

\noindent {\bf Thick Block}\\
Fix constants $D, \epsilon$ and let $\mu = [p,q]$ be  an $\epsilon$-thick
Teichmuller geodesic of length less than $D$. $\mu$ is
$\epsilon$-thick means that for any $x \in \mu$ and any closed
geodesic $\eta$ in the hyperbolic
surface $S^h_x$ over $x$, the length of $\eta$ is greater than
$\epsilon$.
Now let $B^h$ denote the universal curve over $\mu$ reparametrized
such that the length of $\mu$ is covered in unit time. Let $B$ denote
$B^h$ minus a neighborhood of the cusps. Thus $B = S \times [0,1]$
topologically.

A small enough neighborhood of the cusps of $S^h$ is fixed.
 $S^h \times \{ x \}, x \in [0,1]$ is given the hyperbolic structure $S^h_x$
corresponding to the point at distance $x d_{Teich}(p,q)$ from $p$
 along $\mu$ ($d_{Teich}$ denotes Teichmuller metric). A neighborhood
 of the cusps of $S^h$ having been fixed, we remove the images under
 the Teichmuller map (from $S^h_0$ to $S^h_x$) 
of this neighborhood (having first fixed a neighborhood of the cusps
 of $S^h_0$ as the image under the Teichmuller map from $S^h$).

The resultant manifold $B$ (possibly with boundary) is given the path
metric and is 
called a {\bf thick building block}.

Note that after acting by an element of the mapping class group, we
might as well assume that $\mu$ lies in some given compact region of
Teichmuller space. This is because the marking on $S \times \{ 0 \}$
is not important, but rather its position relative to $S \times \{ 1 \}$
Further, since we shall be constructing models only upto
quasi-isometry, we might as well assume that $S^h \times \{ 0 \}$ and $S^h
\times \{ 1 \}$ {\em lie in the orbit} under the mapping class group
of some fixed base surface. Hence $\mu$ can be further simplified to
be a Teichmuller geodesic joining a pair $(p, q)$ amongst a finite set of
points in the orbit of a fixed hyperbolic surface $S^h$.

\smallskip 
 
\noindent {\bf The Model Manifold}\\
 Note that the boundary of a thin block $B_i$ consists of $S \times \{
0,3 \}$ and the intrinsic path metric on each such $S \times \{ 0 \}$
or
 $S \times \{ 3 \}$ is equivalent to the path metric on $S$.
Also,  the boundary of a thick block $B$ consists of $S \times \{
0,1 \}$, where $S^h_0, S^h_1$ lie in some given bounded region of
Teichmuller space. The intrinsic path metrics on each such $S \times \{ 0 \}$
or
 $S \times \{ 1 \}$ is  the path metric on $S$. 

The model
manifold of {\bf i-bounded geometry} is obtained from  $S \times J$
(where $J$ is a sub-interval of $\Bbb{R}$, 
which may be  semi-infinite or bi-infinite. In the former case, we
choose the usual normalisation $J = [ 0, {\infty })$ ) by first choosing a
  sequence
 of blocks $B_i$ (thick or thin) and corresponding intervals $I_i =
 [0,3]$ or $[0,1]$ according as $B_i$ is thick or thin. The
  metric on $S \times I_i$ is then declared to be that on the  building block
  $B_i$. Thus we have,

\smallskip

\noindent {\bf Definition:} {\it A manifold $M$ homeormorphic to $S \times
  J$, where $ J = [0,
  {\infty })$ or $J = ( - \infty , \infty )$, is said to be a model of {\bf
  i-bounded geometry} if  \\
1) there is a fiber preserving homeomorphism from $M$ to
  $\widetilde{S} \times J$ 
that lifts to  a quasi-isometry of universal covers \\ 
2)  there exists a sequence $I_i$ of intervals (with disjoint
 interiors)
 and blocks $B_i$
where the metric on $S \times I_i$ is the same as
 that on some  building block $B_i$ \\
3)  $\bigcup_i I_i = J$ }

\smallskip

The figure below illustrates schematically what the model looks
like. Filled squares correspond to torii along which hyperbolic Dehn
surgery is performed. The blocks which have no filled squares are the
{\it thick blocks} and those with filled squares are the {\it thin
  blocks}

\begin{center}

\includegraphics[height=4cm]{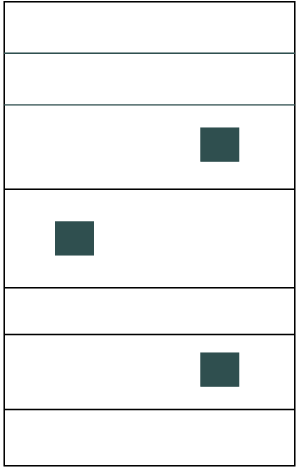}

\smallskip

\underline{Figure 1: {\it Model of i-bounded geometry (schematic)} }

\end{center}

\noindent {\bf Definition:} A manifold $M$ homeormorphic to $S \times
  J$, where $ J = [0,
  {\infty })$ or $J = ( - \infty , \infty )$, is said to have {\bf
  i-bounded geometry} if there exists $K, \epsilon > 0$ such that
 the universal cover $\widetilde{M}$ is $K, \epsilon$ quasi-isometric
  to
a model manifold of i-bounded geometry.

\noindent {\bf The Punctured Torus}\\
In \cite{minsky-torus}, Minsky constructs a model manifold for
arbitrary punctured torus groups that motivates the above
definitions. For him, $S^h$ is the square punctured
torus. $\mathcal{C}$ consists of precisely two shortest curves $a, b$ of
equal length on $S^h$. 
${\mathcal{C}}_i$ is the singleton set $\{ a \}$ for $i$ even and the
set $\{ b \}$ for $i$ odd. The numbers $n$ corresponding to the
surgery coefficients correspond to the number of Dehn twists performed
about the $i$th curve. Thus, we see from Minsky's construction of the
model manifold for
 punctured torus groups that all punctured torus groups give rise to
 manifolds of {\em i-bounded geometry}

\noindent {\bf Alternate Description of i-bounded geometry} \\
We could weaken the definition of thin blocks in models of i-bounded
geometry by requiring that a family $\mathcal{C}$  of disjoint simple
closed curves (rather than a
single simple closed curve) are modified by hyperbolic Dehn
surgery. This gives rise to an equivalent definition.

To see this, add on as many blocks (indexed by $j$)  of $S \times I$
as there are curves 
in $\mathcal{C}$ (this number is bounded in terms of the genus of
$S$). Then isotope Margulis tubes to different levels by a bi-Lipschitz
map away from the tubes. The universal covers of the original $S
\times I$ and the new $S \times \bigcup_j I_j$ are quasi-isometric. 

Hence, it does not multiply examples to allow a family $\mathcal{C}$
rather than a single curve.

\section{Relative Hyperbolicity}

In this section, we shall recall first certain notions of relative
hyperbolicity due to Farb \cite{farb-relhyp}. Using these, we shall
derive certain Lemmas that will be useful in studying the geometry of
the universal covers of building blocks.

\subsection{Electric Geometry}

Given some $\sigma$, we construct a pseudometric, 
  on $S$ by
defining \\
$\bullet$ the length of any path that lies along  $\sigma$ to be zero,\\
$\bullet$  the length of any path $[a,b]$  that misses all
such geodesics in its interior $(a,b)$ to be the hyperbolic length,
and \\
$\bullet$ the length of any other path to be the sum of lengths of
pieces of the above two kinds. \\
This allows us to define distances by taking the infimum of lengths of
  paths joining pairs of points and gives us a path pseudometric,
  which we call the {\bf electric metric}. The electric metric  also
  allows us to define
geodesics. Let us
call $S$ equipped with the  above pseudometric $S_{el}$.

We shall be interested in the universal cover $\widetilde{S_{el}}$ of
$S_{el}$. Paths in $S_{el}$ and $\widetilde{S_{el}}$ will be called
electric paths 
(following Farb \cite{farb-relhyp}). Geodesics and quasigeodesics in 
the electric metric will be called electric geodesics and electric
quasigeodesics respectively.

\smallskip

{\bf Definition:}  $\gamma$ is said to be an electric
$K, \epsilon$-quasigeodesic in $\widetilde{S_{el}}$
{ \bf without backtracking } if
 $\gamma$ is an electric $K$-quasigeodesic in $\widetilde{S_{el}}$ and
 $\gamma$ does not return to any
any lift  $N_\epsilon ( \widetilde{\sigma} ) \subset
\widetilde{S_{el}}$ of
$N_\epsilon ( \sigma )$ after leaving it.

\smallskip

A hyperbolic geodesic $\lambda$ may follow a lift $\widetilde{\sigma}$
for a long time {\it without/ after/ before/ before and after }
intersecting it. This is why in the definition of quasigeodesics
without backtracking, we take $N_\epsilon ( \widetilde{\sigma })$
rather than $\widetilde{\sigma}$ itself.

A similar definition  can be given in the case of manifolds with
cusps. Here electrocuted sets correspond to horodisks (lifts of
cusps). More generally, we can consider $X$ to be a convex subset of
${\Bbb{H}}^n$ and $\mathcal{H}$ to be a collection of uniformly
separated horoballs in $X$ based on  points of $\partial X$ (i.e. they
are the intersection with $X$ of certain horoballs in ${\Bbb{H}}^n$ whose
boundary point lies in $\partial X$). We present below two basic Lemmas due to Farb \cite{farb-relhyp}
in the general setup of hyperbolic metric spaces. Their specializations for  $\widetilde{S_{el}}$
are also indicated.

Let $X$ be a hyperbolic metric space  and $\mathcal{H}$ a collection 
of
{\em (uniformly) $ C$-quasiconvex uniformly separated subsets}, i.e.
there exists $D > 0$ such that for $H_1, H_2 \in \mathcal{H}$, $d_X (H_1,
H_2) \geq D$. In this situation $X$ is hyperbolic relative to the
collection $\mathcal{H}$ (see \cite{bowditch-relhyp}).

\smallskip

{\bf Definition:} A collection $\mathcal{H}$ of uniformly
$C$-quasiconvex sets in a $\delta$-hyperbolic metric space $X$
is said to be {\bf mutually D-cobounded} if 
 for all $H_i, H_j \in \mathcal{H}$, $\pi_i
(H_j)$ has diameter less than $D$, where $\pi_i$ denotes a nearest
point projection of $X$ onto $H_i$. A collection is {\bf mutually
  cobounded } if it is mutually D-cobounded for some $D$. 

\smallskip

{\em Mutual coboundedness} was proven for horoballs by Farb in Lemma 4.7 of
\cite{farb-relhyp} and by Bowditch in stating that the projection of
the link of a vertex onto another \cite{bowditch-relhyp} has bounded
diameter in the link. However, the comparability of 
intersection patterns in this context needs to be stated a bit more
carefully. We give the general version of Farb's theorem below and
refer to \cite{farb-relhyp} \cite{bowditch-relhyp} and Klarreich \cite{klarreich} for proofs.

\begin{lemma} (See Lemma 4.5 and Proposition 4.6 of \cite{farb-relhyp})
Given $\delta , C, D$ there exists $\Delta$ such that
if $X$ is a $\delta$-hyperbolic metric space with a collection
$\mathcal{H}$ of $C$-quasiconvex $D$-separated sets.
then,

\begin{enumerate}
\item {\it Electric quasi-geodesics electrically track hyperbolic
  geodesics:} Given $P > 0$, there exists $K > 0$ with the following
  property: Let $\beta$ be any electric $P$-quasigeodesic from $x$ to
  $y$, and let $\gamma$ be the hyperbolic geodesic from $x$ to $y$. 
Then $\beta \subset N_K^e ( \gamma )$. \\
\item $\gamma$ lies in a {\em hyperbolic} $K$-neighborhood of $N_0 ( \beta
  )$, where $N_0 ( \beta )$ denotes the zero neighborhood of $\beta$
  in the {\em electric metric}. \\
\item {\it Hyperbolicity:} 
  $X$  is $\Delta$-hyperbolic. \\
\end{enumerate}
\label{farb1A}
\end{lemma}

We shall have need to use Lemma \ref{farb1A} in the special
 case that $X = \widetilde{S}$ 
and where the electric metric on  $\widetilde{S_{el}}$ is obtained as at the beginning of this subsection.

\begin{lemma} 
1) Given $P > 0$, there exists $K > 0$ with the following
  property: For some $\widetilde{S_i}$, 
let $\beta$ be any electric $P$-quasigeodesic without backtracking
from $x$ to
  $y$, and let $\gamma$ be the hyperbolic geodesic from $x$ to $y$. 
Then $\beta \subset N_K^e ( \gamma )$. \\
2) There exists $\delta$ such that each
  $\widetilde{S_{el}}$  is $\delta$-hyperbolic, independent of the
  curve $\sigma \in \mathcal{C}$ whose lifts are electrocuted. \\
\label{farb1}
\end{lemma}

We shall need to give a general definition of geodesics and
quasigeodesics without backtracking. \\

{\bf Definitions:} Given a collection $\mathcal{H}$
of $C$-quasiconvex, $D$-separated sets and a number $\epsilon$ we
shall say that a geodesic (resp. quasigeodesic) $\gamma$ is a geodesic
(resp. quasigeodesic) {\bf without backtracking} with respect to
$\epsilon$ neighborhoods if $\gamma$ does not return to $N_\epsilon
(H)$ after leaving it, for any $H \in \mathcal{H}$. 
A geodesic (resp. quasigeodesic) $\gamma$ is a geodesic
(resp. quasigeodesic) {\bf without backtracking} if it is a geodesic
(resp. quasigeodesic) without backtracking with respect to
$\epsilon$ neighborhoods for some $\epsilon \geq 0$.

\smallskip

{\bf Note:} For strictly convex sets, $\epsilon = 0$ suffices, whereas
for convex sets any $\epsilon > 0$ is enough.

\smallskip

Item (2) in  Lemma \ref{farb1A} is due to Klarreich \cite{klarreich},
where the proof is given for $\beta$ an electric geodesic, but the
same proof goes through for electric quasigeodesics. 

\smallskip

{\bf Note:} For Lemma \ref{farb1A}, the hypothesis is that $\mathcal{H}$
  consists of uniformly quasiconvex, mutually separated sets.
Mutual coboundedness has not yet been used. We introduce {\em
  co-boundedness} in the next lemma.

\begin{lemma} 
Suppose $X$ is a $\delta$-hyperbolic metric space with a collection
$\mathcal{H}$ of $C$-quasiconvex $K$-separated $D$-mutually cobounded
subsets. There exists $\epsilon_0 = \epsilon_0 (C, K, D, \delta )$ such that
the following holds:

Let $\beta$ 
  be an electric $P$-quasigeodesic without backtracking
and $\gamma$ a hyperbolic geodesic,
  both joining $x, y$. Then, given $\epsilon \geq \epsilon_0$
 there exists $D = D(P, \epsilon )$ such that \\
1) {\it Similar Intersection Patterns 1:}  if
  precisely one of $\{ \beta , \gamma \}$ meets an
  $\epsilon$-neighborhood $N_\epsilon (H_1)$
of an electrocuted quasiconvex set
  $H_1 \in \mathcal{H}$, then the length (measured in the intrinsic path-metric
  on  $N_\epsilon (H_1)$ ) from the entry point
  to the 
  exit point is at most $D$. \\
2) {\it Similar Intersection Patterns 2:}  if
 both $\{ \beta , \gamma \}$ meet some  $N_\epsilon (H_1)$
 then the length (measured in the intrinsic path-metric
  on  $N_\epsilon (H_1)$ ) from the entry point of
 $\beta$ to that of $\gamma$ is at most $D$; similarly for exit points. \\
\label{farb2A}
\end{lemma}

We summarise the two Lemmas \ref{farb1A} and \ref{farb2A}
in forms that we shall use:

\noindent $\bullet$ If $X$ is a hyperbolic metric space and
$\mathcal{H}$ a collection of uniformly quasiconvex separated subsets,
then $X$ is hyperbolic relative to the collection $\mathcal{H}$. \\

\noindent $\bullet$ If $X$ is a hyperbolic metric space and
$\mathcal{H}$ a collection of uniformly quasiconvex mutually cobounded
separated subsets,
then $X$ is hyperbolic relative to the collection $\mathcal{H}$ and
satisfies {\em Bounded Penetration}, i.e. hyperbolic geodesics and
electric quasigeodesics have similar intersection patterns in the
sense of Lemma \ref{farb2A}. \\

The relevance of co-boundedness comes from the following Lemma which
is essentially due to Farb \cite{farb-relhyp}.

\begin{lemma}
Let $M^h$ be a hyperbolic manifold of {\em i-bounded geometry}, with
Margulis tubes $T_i \in \mathcal{T}$ and horoballs $H_j \in
\mathcal{H}$. Then the lifts $\widetilde{T_i}$ and $\widetilde{H_j}$
are mutually co-bounded.
\label{coboundedHor&T}
\end{lemma}

The proof given in \cite{farb-relhyp} is for a collection of separated
horospheres, but the same proof works for neighborhoods of geodesics
and horospheres as well.

A closely related  theorem was proved by  McMullen
(Theorem 8.1 of \cite{ctm-locconn}).

As usual, $N_R (Z)$ will denote the $R$-neighborhood of the set $Z$. \\
Let $\mathcal{H}$ be a locally finite collection of horoballs in a convex
subset $X$ of ${\Bbb{H}}^n$ 
(where the intersection of a horoball, which meets $\partial X$ in a point, 
 with $X$ is
called a horoball in $X$).

\smallskip

{\bf Definition:} The $\epsilon$-neighborhood of a bi-infinite
geodesic in ${\Bbb{H}}^n$ will be called a {\bf thickened geodesic}. 

\smallskip

\begin{theorem} \cite{ctm-locconn}
Let $\gamma: I \rightarrow X \setminus \bigcup \mathcal{H}$ be an ambient
$K$-quasigeodesic (for $X$ a convex subset  of ${\Bbb{H}}^n$) and let
$\mathcal{H}$  denote a uniformly separated
collection of horoballs and thickened geodesics.
Let $\eta$ be the hyperbolic geodesic with the same endpoints as
$\gamma$. Let $\mathcal{H}({\eta})$  
be the union of all the horoballs and thickened geodesics
 in $\mathcal{H}$ meeting $\eta$. Then
$\eta\cup\mathcal{H}{({\eta})}$ is (uniformly) quasiconvex and $\gamma
(I) \subset  
B_R (\eta \cup \mathcal{H} ({\eta}))$, where $R$ depends only on
$K$. 
\label{ctm}
\end{theorem}

As in Lemmas \ref{farb1A} and \ref{farb2A}, this theorem goes through
for {\em mutually cobounded separated uniformly quasiconvex sets }
$\bf{H}$.

A special kind of {\em geodesic without backtracking} will be
necessary for universal covers of surfaces with some electric metric.

\smallskip

Let $\lambda_e$ be an electric geodesic in some $(\widetilde{S} , d_e)$
for $\widetilde{S}$ equipped with some electric metric obtained by
electrocuting a collection of {\em mutually cobounded separated}
geodesics. Then, each segment of $\lambda_e$ between electrocuted
geodesics is perpendicular to the electrocuted geodesics that it
starts and ends at. We shall refer to these segments of $\lambda_e$ as
{\bf complementary segments} because they lie in the complement of the
electrocuted geodesics. Let $a_\eta , b_\eta$ be the points at which
$\lambda_e$ enters and leaves the electrocuted (bi-infinite) geodesic
$\eta$. Let $[a,b]_\eta$ denote the geodesic segment contained in
$\eta$ joining $a, b$. Segments like $[a,b]_\eta$ shall be referred to
as {\bf interpolating segments}. The union of the {\em complementary
  segments}
along with the {\em interpolating segments} gives rise to a preferred
representative of geodesics joining the end-points of $\lambda_e$; in
fact it is the unique quasigeodesic without backtracking whose
underlying set represents an electric geodesic joining the end-points
of $\lambda_e$. Such a representative of the class of $\lambda_e$
shall be called the {\bf canonical representative} of
$\lambda_e$. Further, the underlying set of the canonical
representative in the {\em hyperbolic metric} shall be called the {\bf
  electro-ambient representative} $\lambda_q$ of $\lambda_e$. 
Since $\lambda_q$ will turn out to be a hyperbolic quasigeodesic, we
shall
also call it an {\bf electro-ambient quasigeodesic}. See Figure 2
below:

\bigskip

\begin{center}

\includegraphics[height=4cm]{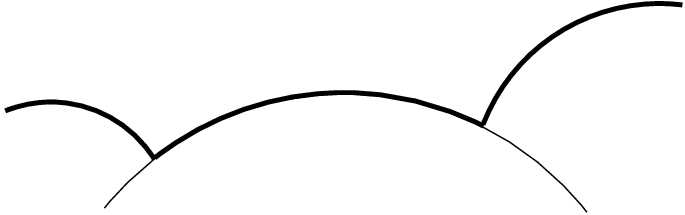}

\underline{Figure 2:{\it Electro-ambient quasigeodesic} }

\end{center}

\smallskip

Now, let $\lambda_h$ denote the hyperbolic geodesic joining the
end-points of $\lambda_e$. By Lemma \ref{farb2A}, $\lambda_h$ and
$\lambda_e$, and hence $\lambda_h$ and $\lambda_q$
have similar intersection patterns with $N_\epsilon (\eta )$ for
electrocuted geodesics $\eta $. Also, $\lambda_h$ and $\lambda_q$
track each other off $N_\epsilon (\eta )$. Further, each {\em
  interpolating segment} of $\lambda_q$ being a {\em hyperbolic} geodesic, it
follows (from the `$K$-fellow-traveller' property of hyperbolic geodesics
starting and ending near each other)
that each interpolating segment of $\lambda_q$
lies within a $(K + 2 \epsilon )$
neighborhood of $\lambda_h$. Again, since each segment of $\lambda_q$
that does not meet an electrocuted geodesic that $\lambda_h$ meets is 
of uniformly bounded (by $C$ say) length, we have finally that
$\lambda_q$ lies within a $(K + C + 2 \epsilon)$ neighborhood of
$\lambda_h$. Finally, since $\lambda_q$ is an electro-ambient
representative, it does not backtrack. Hence we have the following:

\begin{lemma}
There exists $(K, \epsilon)$ such that each electro-ambient
representative of an electric geodesic is a $(K,
\epsilon)$ hyperbolic quasigeodesic.
\label{ea}
\end{lemma}

In the above form, {\em electro-ambient quasigeodesics} are considered
only 
in the context of surfaces and closed geodesics on them. This can be
generalised considerably. Let $X$ be a $\delta$-hyperbolic metric
space, and $\mathcal{H}$ a family of $C$-quasiconvex, $D$-separated,
$k$-cobounded collection of subsets. Then by Lemma \ref{farb1A},
$X_{el}$ obtained by electrocuting the subsets in $\mathcal{H}$ is
a $\Delta = \Delta ( \delta , C, D)$ -hyperbolic metric space. Now,
let $\alpha = [a,b]$ be a hyperbolic geodesic in $X$ and $\beta $ be
an electric 
$P$-quasigeodesic without backtracking
 joining $a, b$. Replace each maximal subsegment
(with end-points $p, q$, say) of
$\beta$ lying within some $H \in \mathcal{H}$
by a hyperbolic {\em interpolating geodesic} $[p,q]$. The resulting
{\bf connected}
path $\beta_q$ is called an {\em electro-ambient quasigeodesic} in
$X$.
The following Lemmas justify the terminology:

\begin{lemma}
Given $\delta, D, C, k, P$ as above, 
there exists $C_3$ such that the following holds: \\
Let $\alpha$, $\beta$ be as above. 
Then $\alpha$ lies in a  $C_3$ neighborhood of $\beta_q$
\label{ea-genl}
\end{lemma}

{\bf Proof:} The proof idea is similar to that for surfaces
and geodesics. 

By Lemma \ref{farb1A}, item (2), there exists $C_0$ such that $\alpha$
lies in a  (hyperbolic) $C_0$-neighborhood of $N_0 ( \beta_q
)$. Further, by {\em bounded penetration} following from
co-boundedness, there exists $C_1$ such that if some interpolating
geodesic $[p,q]$ in $H$ is of length greater than $C_1$, then there
exist $p_1, q_1 \in H \cap \alpha$ such that \\
$d(p,p_1) \leq C_1$ \\
$d(q, q_1) \leq C_1$ \\
$d(p_1, q_1) $ is maximal over all pairs $u, v \in H \cap \alpha$ \\

Hence, by the fellow traveller property, there exists $C_2$ such that
the hyperbolic  geodesic $[p_1, q_1] \subset \alpha$ lies in a
$C_2$-neighborhood of $[p,q]$ and hence $\beta_q$.

Now, if $x \in \alpha$, $x $ lies in a $C_0$ neighborhood of $N_0 (
\beta_q )$. Let $y \in N_0 ( \beta_q )$ be the point nearest to
$x$. If $y$ lies on $\beta - \mathcal{H}$, then $d(x, \beta_q ) \leq
C_0$. Else, $y$ lies on some $H$. Two cases arise:

{\em Case 1:} $\beta$ and hence $\beta_q$ do not penetrate $H$ for
more than $C_1$. In this case, there exists 
 $y \in \beta - \mathcal{H}$, such that $d(x,y) \leq C_0 + C_1$. 

{\em Case 2:}  $\beta$ and hence $\beta_q$ do  penetrate $H$ for
more than $C_1$ and therefore an interpolating geodesic $[p,q]$ of length
greater than $C_1$ exists. Hence there exists a maximal subsegment of
$\alpha$ within a $C_2$ neighborhood of $[p,q]$. From this it follows
easily that $x$ lies in a $C_2$ neighborhood of $\beta_q$.

Thus $\alpha$ lies in a
(uniformly) bounded $C_3$-neighborhood of $\beta_q$. (Here, $C_3 = C_0
+ C_1 + C_2$ suffices.) $\Box$

\smallskip

In
fact, more is true. $\beta_q$ is a hyperbolic quasigeodesic. But to
see this needs a bit more work.
For the sake of concreteness, and to simplify the exposition, we
assume that $X$ is a complete simply connected manifold of pinched
negative curvature. Let $\pi_\alpha$ denote the nearest 
point retraction onto $\alpha$. Since $\beta_q$ is connected, joins the
end-points of $\alpha$ and
$\pi_\alpha$ is continuous, $\pi_\alpha ( \beta_q ) = \alpha$.

{\bf Claim:} There exists $D > 0$ such that
any two points $u, v $ with $\pi_\alpha (u) = \pi_\alpha
(v) = w$ satisfy $d(u,v) \leq D$.

{\bf Proof of Claim:} The  loop that goes from $w$ to $u$ by a
hyperbolic path of length less than $C_3$ (from Lemma 
\ref{ea-genl} ), then from $u$ to $v$ along $\beta_q$ and then back to
$w$ by a hyperbolic path of less than $C_3$ is a $C_4$-quasigeodesic
(for some uniform $C_4$). It can also be converted into a path without
backtracking. Clearly, the geodesic joining the end-points of the path
(a loop)
has length zero. By Lemma \ref{farb2A} the path penetrates each $H \in
\mathcal{H}$ by a uniformly bounded amount $C_5$. Hence, there exists
a uniform constant 
$C_6$ depending on $C_5$ and $D$, such that the loop has (hyperbolic)
length less than $C_6$. The Claim follows. $\Box$

In the same way, it follows that given $D_0$ there exists $D_1$, such
that if $d( \pi_\alpha (u), \pi_\alpha (v)) \leq D_0$ then $d(u,v)
\leq D_1$.

From the Claim above, it follows that $\beta_q$ also must lie in a
bounded neighborhood of $\alpha$ (else there will have to be long
detours along $\beta_q$ starting and ending at a distance less than
$2C_3$ from each other). Further, $\beta_q$ cannot have long pieces
starting and ending close to each other for the same reason. 
Thus $\beta_q$ lying in a bounded neighborhood of $\alpha$ must
`progress'. In other words $\beta_q$ must be a hyperbolic quasigeodesic.
We state this formally below:

\begin{lemma}
There exist $K, \epsilon$ depending on 
 $\delta, D, C, k, P$, such that $\beta_q$ is a $(K, \epsilon
 )$-quasigeodesic. 
\label{ea2}
\end{lemma}

In our proof of Lemma \ref{ea-genl}, we have used the hypothesis that
the collection $\mathcal{H}$ of qc sets is a mutually cobounded
collection. However, this hypothesis can be relaxed. The proof is
exactly the same as Klarreich's  proof of Proposition 4.3 of
\cite{klarreich}, which has been stated here as Item (2) in Lemma
\ref{farb1A} above. We state this below and refer to Proposition 4.3
of \cite{klarreich} for the relevant details:

\begin{lemma}
Given $\delta$, $C, D, P$ there exists $C_3$ such that the following
holds: \\
Let $(X,d)$ be a $\delta$-hyperbolic metric space and $\mathcal{H}$ a
family of $C$-quasiconvex, $D$-separated collection of quasiconvex
subsets. Let $(X,d_e)$ denote the electric space obtained by
electrocuting elements of $\mathcal{H}$.  Then, if $\alpha , \beta_q$
denote respectively a hyperbolic geodesic and an electro-ambient
$P$-quasigeodesic with the same end-points, then $\alpha$ lies in a
(hyperbolic) 
$C_3$ neighborhood of $\beta_q$.
\label{ea-strong}
\end{lemma}

{\bf Note:} The above Lemma generalises Klarreich's Property (2) in
Lemma \ref{farb1A}) by replacing $N_0 ( \beta )$ with $\beta_q$. The
former set can be quite large, but $\beta_q$ is much smaller,
containing only one geodesic segment in $H$ rather than all of $H$. It
is the introduction of the notion of electro-ambient quasigeodesic
that makes for this generalisation. However, Lemma \ref{ea2} is false
in this generality. The idea is that two elements of $\mathcal{H}$
might have geodesics that are parallel  (i.e. close to each other)
for their entire length. Then an electro-ambient quasigeodesic might
look like two adjacent edges of a thin triangle. This is precluded in
Lemma \ref{ea2} by the hypothesis of coboundedness. We shall not be
needing this stronger  Lemma \ref{ea-strong} in this paper and it is
included here for completeness.

\smallskip

Another Lemma which we shall be using follows from the
proof of the Claim in the proof of Lemma \ref{ea2} above.

\begin{lemma}
Given $D_0$ there exists $D_1$ such that 
if $\alpha$ be a loop without backtracking in $X_{el}$ with electric length
less than $D_0$ and further, if $\alpha \cap H$ is a geodesic
for each $H \in \mathcal{H}$, then the hyperbolic length of $\alpha$
is less than $D_1$.
\label{loop-bdd}
\end{lemma}

\subsection{Dehn twists are electric isometries}

Let $S_i$ be a  surface whose path-pseudometric is obtained
from a (fixed) hyperbolic metric by electrocuting the geodesic
$\sigma_i$ in
${\mathcal{C}}$. We can think 
of the Dehn twists as supported in the $\epsilon$-neighborhood
$N_{\epsilon} ( \sigma_i )$ and
that these neighborhoods have been given the zero-metric. Denote the
resultant electric metric on $S_i$ by $\rho_i$

We want to show that any power of a Dehn twist about
$\sigma_i$
induces an isometry of the surface $S_i$ equipped with
$\rho_i$. Consider any two points $x, y \in S$. Let $\alpha$ be any
path in general position with respect to $\sigma$
joining $x, y$.  Look at the action of Dehn twist {\em tw} about
some curve $\sigma_i \in {\mathcal{C}}$ on $\alpha$. Let $\alpha$ meet
$\sigma$ in $p_1, \cdots p_k$. Let $tw( \alpha )$ be the path obtained
from $\alpha$ by keeping it unchanged off $\sigma$ and for each
intersection
point $p_i$, we compose $\alpha$ with a path lying on  $\sigma$
starting and ending at $p_i$ and traversing $\sigma$ once in the
direction of the Dehn twist. Since the restriction $\rho_i|_\sigma =
0$, $\alpha$ and $tw ( \alpha )$ have the same length. Hence, the
length
of the shortest path (geodesic) 
in the homotopy class (rel. end-points) of $tw (
\alpha )$ is less than or equal to the length of the geodesic
reperesentative of the class of $\alpha$. 

Again, let $\beta$ be any path in the homotopy class of $tw ( \alpha
)$. Then by acting by the reverse Dehn twist $tw^{-1}$ about $\sigma$,
we find by an identical argument that the geodesic representative of 
of the homotopy class of $\alpha $, which is the same as that of
$tw^{-1} ( \beta )$ has length less than or equal to the length of the
geodesic representative of $\beta$. 

Since  $\beta$ and $tw ( \alpha )$
are homotopic rel. endpoints, we conclude that $\alpha$ and $tw (
\alpha )$ have geodesic reperesentatives of the same length. 

This proves \\

\begin{lemma}
Let $tw_i^n$ denote a power of a Dehn twist about the curve
$\sigma_i \in {\mathcal{C}}$ and $\rho_i$ denote the electric metric
on $S_i$. Then $tw_i^n $ induces an isometry of $(S_i, \rho_i)$. In
particular, we may arrange $tw_i^n$ to take geodesics to geodesics.
\label{tw-isom1}
\end{lemma}

The last statement in Lemma \ref{tw-isom1} has been put because
geodesics are not uniquely defined in the usual sense in the electric
metric.
But a preferred path does exist, viz. the path which does not
backtrack (or double back) on any $\sigma \in {\mathcal{C}}_i$,
i.e. restricted to $\sigma$ the path is a geodesic in the ordinary
sense.

Everything in the above can be lifted to the universal cover
$\widetilde{S_i}$. We  let $\widetilde{tw}$ denote the lift of $tw$ to
$\widetilde{S_i}$.
This gives \\

\begin{lemma}
Let $\widetilde{tw_i^n}$ denote a lift of $tw_i^n$ to $\widetilde{S_i}$.
Let $\widetilde{\rho_i}$ denote the lifted electric metric
on $\widetilde{S_i}$. Then $\widetilde{tw_i^n} $ induces an isometry of
$({\widetilde{S_i}}, \widetilde{\rho_i})$.  In
particular, we may arrange $\widetilde{tw_i^n}$ to take geodesics to geodesics.
\label{tw-isom2}
\end{lemma}

\subsection{Nearest-point Projections}

We need the following basic lemmas from \cite{mitra-trees}.
The following Lemma  says nearest point projections in a $\delta$-hyperbolic
metric space do not increase distances much.

\begin{lemma}
(Lemma 3.1 of \cite{mitra-trees})
Let $(Y,d)$ be a $\delta$-hyperbolic metric space
 and  let $\mu\subset{Y}$ be a $C$-quasiconvex subset, 
e.g. a geodesic segment.
Let ${\pi}:Y\rightarrow\mu$ map $y\in{Y}$ to a point on
$\mu$ nearest to $y$. Then $d{(\pi{(x)},\pi{(y)})}\leq{C_3}d{(x,y)}$ for
all $x, y\in{Y}$ where $C_3$ depends only on $\delta, C$.
\label{easyprojnlemma}
\end{lemma}

The next lemma  says that quasi-isometries and nearest-point projections on
hyperbolic metric spaces `almost commute'.

\begin{lemma}
(Lemma 3.5 of \cite{mitra-trees})Suppose $(Y_1,d_1)$ and $(Y_2,d_2)$
are $\delta$-hyperbolic.
Let $\mu_1$ be some geodesic segment in $Y_1$ joining $a, b$ and let $p$
be any vertex of $Y_1$. Also let $q$ be a vertex on $\mu_1$ such that
${d_1}(p,q)\leq{d_2}(p,x)$ for $x\in\mu_1$. 
Let $\phi$ be a $(K,{\epsilon})$ - quasiisometric embedding
 from $Y_1$ to $Y_2$.
Let $\mu_2$ be a geodesic segment 
in $Y_2$ joining ${\phi}(a)$ to ${\phi}(b)$ . Let
$r$ be a point on $\mu_2$ such that ${d_2}({\phi}(p),r)\leq{d_2}({\phi(p)},x)$ for $x\in\mu_2$.
Then ${d_2}(r,{\phi}(q))\leq{C_4}$ for some constant $C_4$ depending   only on
$K, \epsilon $ and $\delta$. 
\label{almost-commute}
\end{lemma}

{\bf Sketch of Proof:} (See \cite{mitra-trees} for details.) 
$[p,q]\cup\mu_1$ is called a {\bf tripod}. Then $[p,q]\cup [q,b]$,
$[p,q]\cup [q,a]$ and $[a,b]$ are all quasigeodesics. Hence after
acting by $\phi$ they map to quasigeodesics. In particular, $\phi (q)$
must lie close to the image under $\phi$ of each of 
$[p,q]\cup [q,b]$,
$[p,q]\cup [q,a]$ and $[a,b]$. Hence it must lie close to each of
$[\phi (a), \phi (b)]$, $[\phi (a), \phi (p)]$ and $[\phi (b), \phi
  (p)]$. Again, if $\phi (a), \phi (b), \phi (p), z$ form the four
points of a {\it tripod} (where $z$ is a nearest point projection of
$\phi (p)$ onto the geodesic joining $\phi (a), \phi (b)$), then $z$
too must lie close to each of 
$[\phi (a), \phi (b)]$, $[\phi (a), \phi (p)]$ and $[\phi (b), \phi
  (p)]$.

The result follows by thinness of hyperbolic triangles. $\Box$

\medskip

For our purposes we shall need the above Lemma for quasi-isometries
from $\widetilde{S_a}$ to  $\widetilde{S_b}$ for two different
hyperbolic structures on the same surface. We shall also need it for
the electrocuted surfaces obtained in Lemma \ref{farb1}.

Yet another property that we shall require for nearest point
projections is that nearest point projections in the electric metric
and in the hyperbolic metric almost agree. Let $\widetilde{S} = Y$ be
the universal cover of a surface with the hyperbolic metric minus a
neighborhood of cusps. Equip $Y$ with the path metric $d$
as usual. Then
$Y$ is either the hyperbolic plane (if $S$ has no cusps) or else is
quasi-isometric to a tree (the Cayley graph of a free group). Let
$\sigma$ be a closed geodesic on $S$. Let $d_{e}$ denote the electric
metric on $Y$ obtained by electrocuting the lifts of $\sigma$. Now,
let $\mu = [a,b]$ be a hyperbolic geodesic on $(Y,d)$ and let
$\mu_q$ denote the electro-ambient quasigeodesic 
joining $a, b$. Let $\pi$ denote the nearest point projection in
$(Y,d)$. Tentatively, let $\pi_e$ denote the nearest point projection in
$(Y,d_e)$. Note that $\pi_e$ is not well-defined. It is defined upto a
bounded amount of discrepancy in the electric metric $d_e$. But we
would like to make $\pi_e$ well-defined upto a bounded amount of
discrepancy in the hyperbolic metric $d$.

\smallskip

{\bf Definition:} Let $y \in Y$ and $\mu_q$ be an electro-ambient
representative of an electric geodesic $\mu_e$ in
$(Y,d_e)$. Then $\pi_e(y) = z \in \mu_q$ if the ordered pair $\{
d_e(y,\pi_e(y)), d(y, \pi_e(y) ) \}$ is minimised at $z$. 

Note that this gives us a definition of $\pi_e$ which is ambiguous by
a finite amount of discrepancy not only in the electric metric but
also
in the hyperbolic metric.

\begin{lemma} There exists $C > 0$ such that the following holds.
Let $\mu$ be a hyperbolic geodesic joining $a, b$. Let
  $\mu_e$ be the canonical representative of the electric geodesic
  joining $a, b$. Also let $\mu_q$ be the electro-ambient
  representative of $\mu_e$. Let $\pi_h$ denote the nearest point
  projection of ${\Bbb{H}}^2$ onto $\mu$. 
$d(\pi_h(y) , \pi_e(y))$ is uniformly bounded.
\label{hyp=elproj}
\end{lemma}

{\bf Proof:} 
The proof is similar to that of Lemma \ref{almost-commute}, i.e. Lemma
3.5 of \cite{mitra-trees}. 

$[u, v]_h$ and $[u,v]_e$ will denote respectively the hyperbolic
geodesic and the canonical representative of the electric geodesic
joining $u, v$

$[y,\pi_e(y)] \cup [\pi_e(y),a]$ is an electric quasigeodesic without
backtracking. Hence as in the proof of Lemma \ref{almost-commute},  
$[y,\pi_e(y)] \cup [\pi_e(y),a]$ lies in a bounded neighborhood of
$[y,a]_h$. In particular $\pi_e(y)$ lies in a bounded (hyperbolic)
neighborhood of
$[y,a]_h$. By an identical argument $\pi_e(y)$ 
lies in a bounded neighborhood of
$[y,b]_h$. Again, since $\pi_e(y)$ lies on $\mu_e$, therefore by Lemma
\ref{ea}, $\pi_e(y)$ lies in a bounded neighborhood of $\mu$. Hence
there exists $C > 0$ such that
$\pi_e(y) \in N_C([y,a]_h)\cap N_C([y,b]_h) \cap N_C(\mu )$.

Again, 
$[y,\pi_h(y)] \cup [\pi_h(y),a]$ is a hyperbolic quasigeodesic. Hence   
$[y,\pi_h(y)] \cup [\pi_h(y),a]$ lies in a bounded neighborhood of
$[y,a]_h$. In particular $\pi_h(y)$ lies in a bounded (hyperbolic)
neighborhood of
$[y,a]_h$. By an identical argument $\pi_h(y)$ 
lies in a bounded neighborhood of
$[y,b]_h$. Again, since $\pi_h(y)$ lies on $\mu$, therefore, trivially
$\pi_h(y)$ lies in a bounded neighborhood of $\mu$. Hence
there exists $D > 0$ such that
$\pi_h(y) \in N_D([y,a]_h)\cap N_D([y,b]_h) \cap N_D(\mu )$.
Next, by hyperbolicity (thin-triangles)
$N_D([y,a]_h)\cap N_D([y,b]_h) \cap N_D(\mu )$ and
$N_C([y,a]_h)\cap N_C([y,b]_h) \cap N_C(\mu )$ have diameter bounded by
some $D_1$ depending on $D, C$ and choosing $D = C = max(C,D)$, we get 
$d( \pi_h(y) , \pi_e(y)) \leq D_1$. $\Box$.

\section{Universal Covers of Building Blocks and Electric Geometry}

For most of this section 
(except the last subsection) we shall restrict our attention to
closed surfaces and models corresponding to them. Let $S = S^h$ be a closed
surface with some hyperbolic structure. For surfaces with punctures
$S$ will denote $S^h$ minus a neighborhood of cusps. This will call
for some modifications of the exposition, but {\em not} the overall
construction. Hence, for ease of exposition, we postpone dealing with
cusps till the last subsection of this section.

\subsection{Graph Model of Building Blocks}

{\bf Thin Blocks}

\smallskip

Given a geodesic segment $\lambda \subset \widetilde{S}$ and a basic 
{\it thin building block } $B$, let $\lambda = [a,b] \subset
\widetilde{S} \times \{ 0
\}$ be a geodesic segment, where $\widetilde{S} \times \{ 0
\} \widetilde{B}$, and $B$ is obtained from $S \times I$ by hyperbolic
$(1,n)$ Dehn surgery on $N_\epsilon ( \sigma ) \times [1,2]$.

We shall now build a graph model for $\widetilde{B}$ which will be
quasi-isometric to an electrocuted version of the original model,
where lifts of the curves $\sigma \in \mathcal{C}$ 
which correspond to cores of Margulis tubes
are electrocuted.

On $\widetilde{S} \times \{ 0 \} $ and $\widetilde{S} \times \{ 3 \}$
put the hyperbolic metric obtained
from $S = S^h$. On $\widetilde{S} \times \{ 1 \}$ and $\widetilde{S}
\times  \{ 2 \}$ put the electric metric obtained by electrocuting the
lifts of $\sigma$. This constructs $4$ `sheets' of $\widetilde{S}$
comprising the `horizontal skeleton' of the `graph model' of
$\widetilde{B}$. Now for the vertical strands. On each vertical
element of the form $x \times [0,1]$ and $x \times [2,3]$ put the
Euclidean metric. 

The resulting copy of $\widetilde{B}$ will be called the {\bf graph
  model of a thin block}.

Next, let $\phi$ denote the map induced on $\widetilde{S}$ by
$tw^n_\sigma$, the $n$-fold Dehn twist along $\sigma$. Join each $x
\times \{ 1 \} $ to $\phi (x) \times \{ 2 \}$  by a Euclidean segment
of length $1$. 

\smallskip

{\bf Thick Block}

\smallskip

For a thick block $B = \widetilde{S} \times [0,1]$, recall that $B$ is
the universal curve over a `thick' Teichmuller geodesic $\lambda_{Teich} =
[a,b]$
of length less than
some fixed $D > 0$. Each $S \times \{ x \}$ is identified with the
hyperbolic surface over $(a + \frac{x}{b-a} )$ (assuming that the
Teichmuller geodesic is parametrized by arc-length). 

Here $S \times \{ 0 \} $ is identified with the hyperbolic surface
corresponding to $a$, $S \times \{ 1 \}$ is identified with the
hyperbolic surface corresponding to $b$ and each $(x,a)$ is joined to
$(x,b)$ by a segment of length $1$. 

The resulting model of $\widetilde{B}$ is called a {\bf graph model of
  a thick block}.
\smallskip

{\bf Admissible Paths}

\smallskip

Admissible paths consist of the following: \\
1) Horizontal segments along some $\widetilde{S} \times \{ i \}$
  for $ i = \{ 0,1,2,3 \}$ (thin blocks) or $i = \{ 0, 1 \}$ (thick
  blocks).\\
2)  Vertical segments $x \times [0,1]$ or $x \times [2,3]$ for thin
  blocks
or $x \times [0,1]$ for thick blocks.\\
3)  Vertical segments of length $1$ joining $x  \times \{ 1 \}$ to
$\phi (x)  \times \{ 2 \}$ for thin blocks.

\subsection{Construction of Quasiconvex Sets for Building Blocks}

In the next section, we will construct a set $B_\lambda$ containing
$\lambda$ and a retraction $\Pi_\lambda$ of $\widetilde{M}$ onto
it. $\Pi_\lambda$ will have the property that it does not stretch
distances much. This will show that $B_\lambda$ is quasi-isometrically
embedded in $\widetilde{M}$. 

In this subsection, we describe the construction of $B_\lambda$
restricted to a building block $B$.

\smallskip

\noindent {\bf Construction of $B_\lambda (B)$ - Thick Block}\\
Let the thick block be the universal curve over a Teichmuller geodesic
$[\alpha , \beta ]$. Let $S_\alpha$ denote the hyperbolic surface over
$\alpha$ and $S_\beta$ denote the hyperbolic surface over $\beta$.

First, let $\lambda = [a,b]$ be a geodesic segment in
$\widetilde{S}$. Let $\lambda_{B0}$ denote $\lambda \times \{ 0
\}$. 

Next, let  $\phi$ be  the lift of the 'identity' map from 
$\widetilde{S_\alpha}$ to
$\widetilde{S_\beta}$. 
. Let $\Phi$ denote
the induced map on geodesics and let $\Phi ( \lambda )$ denote the
hyperbolic geodesic 
joining $\phi (a), \phi (b)$. Let $\lambda_{B1}$ denote $\Phi (\lambda
) \times \{ 1 \}$.

For the universal cover $\widetilde{B}$ of the thick block $B$, define:

\begin{center}

$B_\lambda (B) = \bigcup_{i=0,1} \lambda_{Bi}$

\end{center}

{\bf Definition:} Each $\widetilde{S} \times i$ for $i = 0, 1$
 will be called a {\bf horizontal sheet} of
$\widetilde{B}$ when $B$ is a thick block.

\smallskip

{\bf Construction of $B_\lambda (B)$ - Thin Block}

\smallskip

First, recall that $\lambda = [a,b]$ is a geodesic segment in
$\widetilde{S}$. Let $\lambda_{B0}$ denote $\lambda \times \{ 0
\}$. 

Next, let $\lambda_{el}$ denote the electric geodesic joining $a,
b$ in the electric pseudo-metric on $\widetilde{S}$ obtained by
electrocuting lifts of $\sigma$. Let $\lambda_{B1}$ denote
$\lambda_{el} \times \{ 1 \}$. 

Third, recall that $\phi$ is the lift of the Dehn twist $tw^n_\sigma$ to
$\widetilde{S}$ equipped with  the electric metric. Let $\Phi$ denote
the induced map on geodesics, i.e. if $\mu = [x,y] \subset (
\widetilde{S} , d_{el} )$, then $\Phi ( \mu ) = [ \phi (x), \phi (y)
]$ is the geodesic joining $\phi (x), \phi (y)$. Let $\lambda_{B2}$
denote $\Phi ( \lambda_{el} ) \times \{ 2 \}$.

Fourthly, let $\Phi ( \lambda )$ denote the hyperbolic geodesic
joining $\phi (a), \phi (b)$. Let $\lambda_{B3}$ denote $\Phi (\lambda
) \times \{ 3 \}$.

For the universal cover $\widetilde{B}$ of the thin block $B$, define:

\begin{center}

$B_\lambda (B) = \bigcup_{i=0,\cdots , 3} \lambda_{Bi}$

\end{center}

{\bf Definition:} Each $\widetilde{S} \times i$ for $i = 0
\cdots 3$ will be called a {\bf horizontal sheet} of
$\widetilde{B}$ when $B$ is a thick block.

\smallskip

{\bf Construction of $\Pi_{\lambda ,B}$ - Thick Block}

\smallskip

On
$\widetilde{S} \times \{ 0 \}$, let $\Pi_{B0}$ denote nearest point
projection
onto $\lambda_{B0}$ in the path metric on $\widetilde{S} \times \{ 0 \}$. 

On
$\widetilde{S} \times \{ 1 \}$, let $\Pi_{B1}$ denote nearest point
projection
onto $\lambda_{B1}$ in the path metric on $\widetilde{S} \times \{ 1 \}$.

For the universal cover $\widetilde{B}$ of the thick block $B$, define:

\begin{center}

$\Pi_{\lambda ,B}(x) = \Pi_{Bi}(x) , x \in \widetilde{S} \times \{ i
  \} , i=0,1$

\end{center}

\smallskip

{\bf Construction of $\Pi_{\lambda ,B}$ - Thin Block}

\smallskip

On 
$\widetilde{S} \times \{ 0 \}$, let $\Pi_{B0}$ denote nearest point
projection onto $\lambda_{B0}$. Here the nearest point projection is
taken in the path metric on $\widetilde{S} \times \{ 0 \}$ which is a
hyperbolic 
metric space.

On $\widetilde{S} \times \{ 1 \}$, let $\Pi_{B1}$ 
denote  the nearest point projection onto $\lambda_{B1}$. Here the
nearest point projection is taken in the sense of the definition
preceding
Lemma \ref{hyp=elproj}, that is minimising the ordered pair $(d_{el},
d_{hyp})$ (where $d_{el}, d_{hyp}$ refer to electric and hyperbolic
metrics respectively.)

On $\widetilde{S} \times \{ 2 \}$, let $\Pi_{B2}$ 
denote  the nearest point projection onto $\lambda_{B2}$. Here, again the
nearest point projection is taken in the sense of the definition
preceding
Lemma \ref{hyp=elproj}.

Again, on 
$\widetilde{S} \times \{ 3 \}$, let $\Pi_{B3}$ denote nearest point
projection onto $\lambda_{B3}$. Here the nearest point projection is
taken in the path metric on $\widetilde{S} \times \{ 3 \}$ 
which is a hyperbolic
metric space.

For the universal cover $\widetilde{B}$ of the thin block $B$, define:

\begin{center}

$\Pi_{\lambda ,B}(x) = \Pi_{Bi}(x) , x \in \widetilde{S} \times \{ i
  \} , i=0,\cdots , 3$

\end{center}

{\bf $\Pi_{\lambda , B}$ is a retract - Thick Block}

The proof for a thick block is exactly as in \cite{mitra-trees}.
The crucial tool is Lemma \ref{almost-commute}. 

\begin{lemma}
There exists $C > 0$ such that the following holds: \\
Let $x, y \in \widetilde{S} \times \{ 0, 1\} \subset \widetilde{B}$
for some thick block $B$. 
Then $d( \Pi_{\lambda , B} (x), \Pi_{\lambda , B} (y)) \leq C d(x,y)$.
\label{retract-thick}
\end{lemma}

{\bf Proof:} It is enough to show this for the two following cases: \\
1) $x, y \in \widetilde{S} \times \{ 0 \} $ OR 
 $x, y \in \widetilde{S} \times \{ 1 \} $. \\
2) $x, y$ are of the form $(p,0)$ and $(\phi (p),1)$ which are
  connected by a vertical segment of length one ( as per the
  construction of the model $B$).

\smallskip

Case 1 above follows from Lemma \ref{easyprojnlemma}, and Case 2 from
the fact that $\phi$ is a uniform quasi-isometry (depending on the
uniform bound on the length of the Teichmuller geodesic over which $B$
is the universal curve) and Lemma \ref{almost-commute} which says that
there exists $C_1 > 0$ such that if $\pi$ be the nearest point
retraction in $\widetilde{S}$ onto $\lambda$ then 
$d(\phi ( \pi (p)), \pi ( \phi (p))) \leq C_1$.
 From this it follows that

\begin{center}

 $d( \Pi_{\lambda , B} ((p,0)), \Pi_{\lambda , B} ((\phi (p), 1)))
 \leq C_1 + 1$

\end{center}

Choosing $C = C_1 + 1$ we are through. $\Box$

{\bf $\Pi_{\lambda , B}$ is a retract - Thin Block}

The two main ingredients in this case are Lemmas \ref{almost-commute}
and \ref{hyp=elproj}. 

\smallskip

{\bf Note:} In the Lemma \ref{retract-thin} below,
there is implicit a constant $n$, the twist coefficient of the Dehn
twist that distinguishes $B$. But the constant $C$ below is
independent of $n$ due to the fact that powers of Dehn twists are {\bf
  uniform quasi-isometries} of the electric metric. In fact this is
the reason why we introduce the electric metric in the first place, so
as to ensure that the techniques of \cite{mitra-trees} and
\cite{mitra-ct} go through here.

\begin{lemma}
There exists $C > 0$ such that the following holds: \\
Let $x, y \in \widetilde{S} \times \{ 0,1,2,3 \} \subset \widetilde{B}$
for some thin block $B$. 
Then $d_e( \Pi_{\lambda , B} (x), \Pi_{\lambda , B} (y)) \leq C d_e(x,y)$.
\label{retract-thin}
\end{lemma}

{\bf Proof:} It is enough to show this for the two following cases: \\
1) $x, y \in \widetilde{S} \times \{ 0 \} $ OR 
 $x, y \in \widetilde{S} \times \{ 1 \} $. \\
2) $x = (p,0)$ and $y = (p,1)$ for some $p$ \\
3) $x, y$ are of the form $(p,1)$ and $(\phi (p),2)$ which are
  connected by a vertical segment of length one ( as per the
  construction of the model $B$) \\
4) $x = (p,2)$ and $y = (p,3)$ for some $p$. \\

\smallskip

{\bf Case 1:} As in Lemma \ref{retract-thick} above, this follows from \ref{easyprojnlemma}

\smallskip

{\bf Case 2 and Case 4:} These follow from Lemma \ref{hyp=elproj} which says
that the hyperbolic and electric projections of $p$ onto the
hyperbolic geodesic $[a,b]$ and the electro-ambient geodesic
$[a,b]_{ea}$ respectively `almost agree'.  If $\pi_h$ and $\pi_e$
denote the hyperbolic and electric projections, then
there exists $C_1 > 0$ such that
$d( \pi_h(p), \pi_e(p)) \leq C_1$.
Hence

\begin{center}

 $d( \Pi_{\lambda , B} ((p,i)), \Pi_{\lambda , B} ((p,i+1)))
 \leq C_1 + 1$, for $i = 0, 2$.

\end{center}

{\bf Case 3:} First, from Lemma \ref{tw-isom2} the (power of the) Dehn twist
 $\phi$ used in the construction of $B$  is a (uniform)
quasi-isometry of $\widetilde{S}$ equipped with the electric
metric. Again, if $\pi$ denotes the nearest point projection in the
 electric metric, then from Lemma \ref{almost-commute}, there exists
 $C_2 > 0$ such that $d_e(\phi ( \pi (p)), \pi ( \phi (p))) \leq C_2$.
 Here $d_e$ denotes the electric metric. From this it follows that

\begin{center}

 $d_e( \Pi_{\lambda , B} ((p,1)), \Pi_{\lambda , B} ((\phi (p), 2)))
 \leq C_1 + 1$

\end{center}

Choosing $C = max (C_1 + 1,C_2 + 1)$ we are through. $\Box$

\subsection{Modifications for Surfaces With Punctures}

We deal with the thin block first.

\smallskip

{\bf Thin Block}

\smallskip

For $S^h$ a hyperbolic surface with punctures, let $S$ denote $S^h$
minus some neighborhood of the cusps. Then the construction of the
model $B$ and hence the graph model of $\widetilde{B}$
for a thin block $B$ goes through  {\em mutatis mutandis} even
with respect to  notation. The construction of the quasi-convex set
and the retraction are modified as follows: \\

\smallskip

$\lambda$ will no longer be a hyperbolic geodesic but rather an
ambient quasigeodesic in $\widetilde{S}$. The construction is taken
from \cite{brahma-pared}. We start with a
hyperbolic geodesic $\lambda^h$ in $S^h$. Fix a neighborhood of the
cusps lifting to an equivariant family of horoballs in the universal
cover ${\Bbb{H}}^2 = \widetilde{S^h}$. 
 Since $\lambda^h$ is a hyperbolic
geodesic in  $\widetilde{S^h}$ there are unique entry and exit
points for each horoball that $\lambda^h$ meets and hence unique
Euclidean geodesics joining them on the corresponding
horosphere. Replacing the segments of $\lambda^h$ lying inside
$Z$-horoballs by the corresponding Euclidean geodesics, we obtain an
ambient quasigeodesic $\lambda$ in $\widetilde{M_0}$ by Theorem
\ref{ctm}. See Figure  below:

\medskip

\begin{center}
\includegraphics{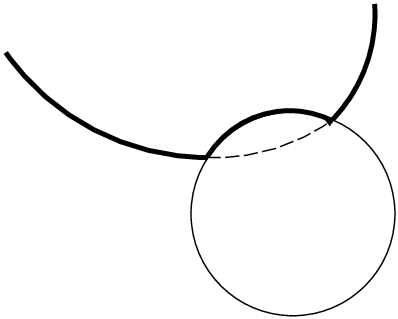}

\smallskip

\underline{Figure 3: Horo-ambient quasigeodesics}

\end{center}

\smallskip

Ambient quasigeodesics obtained by this kind of a construction will be
termed {\bf horo-ambient quasigeodesics} to distinguish them from {\em
  electro-ambient quasigeodesics} defined earlier.
Starting with a horo-ambient quasigeodesic $\lambda \subset
\widetilde{S}$, we can proceed as before to construct $B_{\lambda ,
  B}$, $\Pi_{\lambda , B}$ and prove
 Lemma \ref{retract-thin}
above. In fact the statement and proof of this Lemma goes through unchanged,
with the only pro viso that for punctured surfaces, $S$ and $S^h$ are
not the same and that $\lambda$ is a horo-ambient quasigeodesic in
general.
(Note that if $S$ has no punctures, a horo-ambient quasigeodesic is a
hyperbolic geodesic.)

\smallskip

{\bf Thick Block}

Here $B$ is obtained from the universal curve over a Teichmuller
geodesic by removing a neighborhood of the cusps. Again, $S$ is
obtained from $S^h$ by removing a neighborhood of the cusps. 
That the map $\phi$ is a uniform quasi-isometry is ensured by the
fact that the corresponding copies of $S^h$ are a uniformly bounded
Teichmuller distance from each other, and that if $\phi$ denote a map
between these copies of $S^h$, one can ensure that $\phi$ takes cusps
to cusps. 

The construction of the graph model for $\widetilde{B}$, the
construction of $B_{\lambda , B}$ and $\Pi_{\lambda , B}$ also go
through as before with the pro viso that $\lambda$ is a {\em
  horo-ambient quasigeodesic}. Lemma \ref{retract-thick} goes through
as before.

\section{Construction of Quasiconvex Sets}

\subsection{Construction of $B_\lambda$ and $\Pi_\lambda$}

Given a manifold $M$ of i-bounded geometry, we know that $M$ is
homeomorphic to $S \times J$ for $J = [0, \infty ) $ or $( - { \infty
  }, {\infty })$. By definition of i-bounded geometry, 
there exists a sequence $I_i$ of intervals and blocks $B_i$
where the metric on $S \times I_i$ coincides
  with that on some  building block $B_i$. Denote: \\
$\bullet$ $B_{\mu , B_i} = B_{i \mu }$ \\
$\bullet$ $\Pi_{\mu , B_i} = \Pi_{i \mu }$ \\

Now for a block $B = S \times I$ (thick or thin),  a natural
map $\Phi_B$ may be defined taking
 $ \mu = B_{\mu , B} \cap \widetilde{S} \times \{ 0 \} $ to a
geodesic $B_{\mu , B} \cap \widetilde{S} \times \{ k \} = \Phi_B ( \mu
)$ where $k = 1$ or $3$ according as $B$ is thick or thin. Let the map
$\Phi_{B_i}$ be denoted as $\Phi_i$ for $i \geq 0$.
For $i < 0$ we shall modify this by defining $\Phi_i$ to be the map
that takes 
 $ \mu = B_{\mu , B_i} \cap \widetilde{S} \times \{ k \} $ to a
geodesic $B_{\mu , B_i} \cap \widetilde{S} \times \{ 0 \} = \Phi_i ( \mu
)$ where $k = 1$ or $3$ according as $B$ is thick or thin. 

We start with a reference block $B_0$ and a reference geodesic segment
$\lambda = \lambda_0$ on the `lower surface' $\widetilde{S} \times \{
0 \}$.
Now inductively define: \\
$\bullet$ $\lambda_{i+1}$ = $\Phi_i ( \lambda_i )$ for $i \geq 0$\\
$\bullet$ $\lambda_{i-1}$ = $\Phi_i ( \lambda_i )$ for $i \leq 0$ \\
$\bullet$ $B_{i \lambda }$ = $B_{\lambda_i } ( B_i )$ \\
$\bullet$ $\Pi_{i \lambda }$ = $\Pi_{\lambda_i , B_i}$ \\
$\bullet$ $B_\lambda = \bigcup_i B_{i \lambda }$ \\
$\bullet$ $\Pi_\lambda = \bigcup_i \Pi_{i \lambda }$ \\

Recall that 
each $\widetilde{S} \times i$ for $i = 0
\cdots K$ is called a {\bf horizontal sheet} of
$\widetilde{B}$, where $K = 1$ or $3$ according as $B$ is thick or
thin. We will restrict our attention to the union of the horizontal
sheets $\widetilde{M_H}$ of $\widetilde{M}$ with the induced
metric. 

Clearly, 
$B_\lambda  \subset \widetilde{M_H} \subset \widetilde{M}$, and
$\Pi_\lambda$ is defined from  $\widetilde{M_H}$ to $B_\lambda$. Since
 $\widetilde{M_H}$ is a `coarse net' in $\widetilde{M}$, we will be
able to get all the coarse information we need by restricting
ourselves to  $\widetilde{M_H}$.

By Lemmas \ref{retract-thick} and \ref{retract-thin}, we obtain the
fact that each $\Pi_{i \lambda }$ is a retract. Hence assembling all
these retracts together, we have the following basic theorem:

\begin{theorem}
There exists $C > 0$ such that 
for any geodesic $\lambda = \lambda_0 \subset \widetilde{S} \times \{
0 \} \subset \widetilde{B_0}$, the retraction $\Pi_\lambda :
\widetilde{M_H} \rightarrow B_\lambda $ satisfies: \\

Then $d( \Pi_{\lambda , B} (x), \Pi_{\lambda , B} (y)) \leq C d(x,y) +
C$.
\label{retract}
\end{theorem}

{\bf NOTE 1} For Theorem \ref{retract} above, note that all that we
really
require is that the universal cover $\widetilde{S}$ be a hyperbolic
metric space. There is {\em no restriction on $\widetilde{M_H}$.} In fact,
Theorem \ref{retract} would hold for general stacks of hyperbolic
metric spaces with hyperbolic Dehn surgery performed on blocks. 

{\bf Note 2:} $M_H$ has been given built up out of {\bf graph models
  of thick and thin blocks} and have sheets that are electrocuted.

\subsection{Heights of Blocks}

Instead of considering all the horizontal sheets, we would now like
to consider only the {\bf boundary horizontal sheets}, i.e. for a thick
block we consider $\widetilde{S} \times \{ 0, 1\}$ and for a thin
block we consider $\widetilde{S} \times \{ 0,3 \}$. The union of all 
boundary horizontal sheets will be denoted by $M_{BH}$. 

\smallskip

{\bf Observation 1:} $\widetilde{M_{BH}}$ is a `coarse net' in $\widetilde{M}$
in the {\bf graph model}, but not in the {\bf model of i-bounded geometry}.

In the graph model, any point can be connected by a vertical segment
of length $\leq 2$ to one of the boundary horizontal sheets.

However, in the model of i-bounded geometry, there are points within
Margulis tubes (say for instance, the center of the totally geodesic
disk bounded by a meridian) which are at a distance of the order of $ln
(n_i )$ from the boundary horizontal sheets. Since $n_i$ is arbitrary,
$\widetilde{M_{BH}}$ is no longer a `coarse net' in $\widetilde{M}$
equipped with 
the model of {\it i-bounded geometry}.

{\bf Observation 2:} $\widetilde{M_H}$ is defined only in the {\bf
  graph model}, but not in the model of i-bounded geometry.

{\bf Observation 3:} The electric metric on the model of i-bounded
geometry on $\widetilde{M}$ obtained by
electrocuting all lifts of Margulis tubes is quasi-isometric to the
graph model of $\widetilde{M}$.

This follows from the fact that any lift of a  Margulis tube has
diameter $1$ in the graph model of $\widetilde{M}$.

\smallskip

{\bf Bounded Height of Thick Block}

\smallskip

Let $\mu \subset \widetilde{S} \times \{ 0 \} \widetilde{B_i}$ be a
geodesic in a (thick or thin) block. 
Then  there exists a $(K_i, \epsilon_i )$- quasi-isometry $\psi_i$ ( =
$\phi_i$ for thick blocks)
from $\widetilde{S} \times \{ 0 \}$ to $\widetilde{S} \times \{ 1 \}$
and $\Psi_i$ is the induced map on geodesics. Hence, for any $x \in
\mu$, $\psi_i (x)$ lies within some bounded distance $C_i$ of $\Psi_i
( \mu )$. But $x$ is connected to $\psi_i (x)$ by \\

\smallskip

\noindent {\bf Case 1 - Thick Blocks:} a vertical segment of length $1$  \\
{\bf Case 2 - Thin Blocks:}  the union of \\
1)
two vertical segments of length $1$ between $\widetilde{S}
\times \{ i \}$ and $\widetilde{S} \times \{ i + 1 \}$ for $i = 0, 2$
\\
2) a horizontal segment of length bounded by (some uniform)
  $C^{\prime}$   (cf. Lemma \ref{ea})
connecting $(x,1)$ to a point on the
  electro-ambient geodesic 
  $B_\lambda (B) \cap \widetilde{S} \times \{ 1 \}$ \\
3) a vertical segment of length one in the {\bf graph model}
 connecting $(x,1)$ to $( \phi (x), 2)$. Such a path has to travel
 {\em through the Margulis tube} in the model of {\bf i-bounded
 geometry} and has length less than
$g_0(n_i)$ for some function $g_0: \Bbb{Z}
 \rightarrow \Bbb{N}$, and $n_i$ the twist coefficient.  \\
4) a horizontal segment of length less than $C^{\prime}$ (Lemma \ref{ea})
connecting $(\phi_i (x),3)$ to a point on the
  hyperbolic geodesic 
  $B_\lambda (B) \cap \widetilde{S} \times \{ 3 \}$ \\

\smallskip

Thus $x$ can be connected to a point on $x^{\prime} \in \Psi_i ( \mu
)$ by a path of length less than $g(i) = 2 + 2C^{\prime} + g_0
(n_i)$. Here, in fact $g_0$ is at most linear in $n_i$ but we shall
not need this. Recall that $\lambda_i$ is the geodesic on the lower
horizontal surface of the block $\widetilde{B_i}$. The same can be
done for blocks $\widetilde{B_{i-1}}$ and {\em going down} from
$\lambda_i$ to $\lambda_{i-1}$.
What we have thus shown is:

\begin{lemma}
There exists a function $g: \Bbb{Z} \rightarrow \Bbb{N}$ such that
for any block $B_i$ (resp. $B_{i-1}$), and $x \in \lambda_i$, there exists $x^{\prime} 
\in \lambda_{i+1}$ (resp. $\lambda_{i-1}$) for $i \geq 0$ (resp. $i \leq
0$), satisfying:

\begin{center}

$d(x, x^{\prime}) \leq g(i)$

\end{center}
\label{bddheight}
\end{lemma}

\smallskip

{\bf Modifications for Punctured Surfaces}

For a punctured surface, the above  argument  has to be modified using
some constructions from Lemma 5.1 of \cite{brahma-pared}.

Given $\lambda^h \in \widetilde{S^h}$ we have already indicated
how to construct a horo-ambient quasigeodesic $\lambda$ in
$\widetilde{S}$ (where, recall that $S$ is $S^h$ minus a neighborhood
of cusps). Let {\bf $\lambda_c$} denote the union of the segments of
$\lambda$ that lie along cusps. Let $\lambda_b = \lambda -
\lambda_c$. For punctured surfaces, recall that
$\lambda_i$ is a horo-ambient quasigeodesic on the lower horizontal
surface of $\widetilde{B_i}$. $\lambda_{ic}$ will denote the union of
segments of $\lambda_i$ lying along cusps and $\lambda_{ib}$ will
denote $\lambda_i - \lambda_{ic}$.

Lemma 5.1 of \cite{brahma-pared} says that there exists $C_0$ such
that
for any thick block $B_i$,
and $x \in \lambda_{ib}$ there exists $x^{\prime} \in \lambda_{i+1,
  b}$ such that  $d(x, x^{\prime}) \leq C^{\prime}$. Combining this
with the argument given above for surfaces without punctures, we
conclude,

\begin{lemma}
There exists a function $g: \Bbb{Z} \rightarrow \Bbb{N}$ such that
for any block $B_i$ (resp. $B_{i-1}$), and $x \in \lambda_{ib}$, there exists $x^{\prime} 
\in \lambda_{{i+1},b}$ (resp. $\lambda_{{i-1},b}$) for $i \geq 0$ (resp. $i \leq
0$), satisfying:

\begin{center}

$d(x, x^{\prime}) \leq g(i)$

\end{center}
\label{bddheight-punctured}
\end{lemma}

{\bf Note:} For a surface without punctures, $\lambda$ and $\lambda_b$
coincide.

\section{Cannon-Thurston Maps for Surfaces Without Punctures}

Unless explicitly mentioned otherwise, we shall 
assume till the end of this section that \\
$\bullet$ $S$ is a closed surface. Hence $S^h = S$.\\
$\bullet$ there exists a hyperbolic
manifold $M$ and a homeomorphism from $\widetilde{M}$ to $\widetilde{S} \times
\Bbb{R}$. We identify $\widetilde{M}$ with $\widetilde{S} \times
\Bbb{R}$ via this homeomorphism.\\
$\bullet$  $\widetilde{S} \times \Bbb{R}$ admits a 
quasi-isometry $g$ to a model manifold of {\em i-bounded geometry}\\
$\bullet$
$g$  preserves the fibers over $\Bbb{Z} \subset \Bbb{R}$

\smallskip

{\bf Remarks: 1)} The above assumption is much stronger than what we
need. It suffices to assume that $\widetilde{M}$ is a Gromov-hyperbolic
metric space. Further relaxations on the hypothesis may be
considered while generalising the results of this paper to other
hyperbolic metric spaces.\\
{\bf 2)} We have taken $J$ to be $\Bbb{R}$ here for concreteness. The
other possibility of $J = { \Bbb{R} }^+$ can be treated in exactly the
same way. 

\smallskip

We shall henceforth ignore the quasi-isometry $g$ and think of
$\widetilde{M}$ itself as the universal cover of a model manifold of
{\em i-bounded geometry}.

\subsection{Admissible Paths}

We want to define a collection of {\bf $B_\lambda$-elementary admissible paths}
 lying in a bounded
neighborhood of  $B_\lambda$. $B_\lambda$ is not connected. Hence, it
 does not make much sense to speak of the path-metric on
 $B_\lambda$. To remedy this we introduce a `thickening'
 (cf. \cite{gromov-ai}) of $B_\lambda$ which is path-connected and
 where the
 paths are controlled. A {\bf  $B_\lambda$-admissible path} will be a composition
 of  $B_\lambda$-elementary admissible paths.

Recall
that admissible paths in the graph model of bounded geometry 
consist of the following :\\
1) Horizontal segments along some $\widetilde{S} \times \{ i \}$
  for $ i = \{ 0,1,2,3 \}$ (thin blocks) or $i = \{ 0, 1 \}$ (thick
  blocks).\\
2)  Vertical segments $x \times [0,1]$ or $x \times [2,3]$ for thin
  blocks, where $x \in \widetilde{S}$.\\
3)  Vertical segments of length $1$ joining $x  \times \{ 1 \}$ to
$\phi (x)  \times \{ 2 \}$ for thin blocks.\\
4)  Vertical segments of length $1$ joining $x  \times \{ 0 \}$ to
$\phi (x)  \times \{ 1 \}$ (or $x \times \{ 1 \}$ if $\phi$ is
  regarded as the identity map at the topological level) for thick blocks.

\smallskip

We shall choose a subclass of these admissible paths to define
$B_\lambda$-elementary admissible paths.

\smallskip

\noindent {\bf  $B_\lambda$-elementary admissible paths in the thick block}\\
Let $B = S \times [i,i+1]$ be a thick block, where each $(x,i)$ is
connected by a vertical segment of length $1$ to $( \phi (x),i+1)$.
 Also  $\Phi$ is the map on geodesics induced by
$\phi$. Let $B_\lambda \cap \widetilde{B} = \lambda_i \cup
\lambda_{i+1}$ where $\lambda_i$ lies on $\widetilde{S} \times \{ i
\}$ and $\lambda_{i+1}$ lies on  $\widetilde{S} \times \{ i+1
\}$.  $\pi_j$, for $j = i, i+1$  denote nearest-point projections of
$\widetilde{S} \times \{ j \}$ onto $\lambda_j$. Next, since $\phi$
is a quasi-isometry, there exists $C > 0$ such that for all $(x,i) \in
\lambda_i$, $(\phi(x),i+1)$ lies in a $C$-neighborhood of $\Phi
(\lambda_i ) = \lambda_{i+1}$. The same holds for $\phi^{-1}$ and
points in $\lambda_{i+1}$, where $\phi^{-1}$ denotes the {\em
  quasi-isometric inverse} of $\phi$ from
$\widetilde{S} \times \{ i + 1 \}$ to $\widetilde{S} \times \{ i \}
$. The {\bf
  $B_\lambda$-elementary admissible paths} in $\widetilde{B}$ consist
of the following:\\

1) Horizontal geodesic subsegments of $\lambda_j$,  $j = \{ i, i+ 1 \}$.\\
2)  Vertical segments of length $1$ joining $x  \times \{ 0 \}$ to
$\phi (x)  \times \{ 1 \}$.\\
3)  Vertical segments of length $1$ joining $y  \times \{ 1 \}$ to
$\phi^{-1} (y)  \times \{ 0 \}$.\\
4)  Horizontal geodesic segments lying in a $C$-neighborhood of
  $\lambda_j$, $j = i, i+1$.

\smallskip

{\bf  $B_\lambda$-elementary admissible paths in the thin block}

Let $B = S \times [i,i+3]$ be a thin block, where each $(x,i+1)$ is
connected by a vertical segment of length $1$ to $( \phi (x),i+2)$.
 Also  $\Phi$ is the map on canonical representatives of electric geodesics induced by
$\phi$. Let $B_\lambda \cap \widetilde{B} = \bigcup_{j=i \cdots i+3}
 \lambda_j$ where $\lambda_j$ lies on $\widetilde{S} \times \{ j 
\}$.  $\pi_j$  denotes nearest-point projection of
$\widetilde{S} \times \{ j \}$ onto $\lambda_j$ (in the appropriate
sense - hyperbolic for $j = i, i+3$ and electric for $j = i+1,
i+2$). Next, since $\phi$ 
is an electric isometry, but a hyperbolic quasi-isometry, there exists
$C > 0$ (uniform constant) and $K=K(B)$  such that for all $(x,i) \in
\lambda_i$, $(\phi(x),i+1)$ lies in an (electric) $C$-neighborhood and
a hyperbolic $K$-neighborhood of $\Phi
(\lambda_{i+1} ) = \lambda_{i+2}$. The same holds for $\phi^{-1}$ and
points in $\lambda_{i+2}$, where $\phi^{-1}$ denotes the {\em
  quasi-isometric inverse} of $\phi$ from
$\widetilde{S} \times \{ i + 2 \}$ to $\widetilde{S} \times \{ i + 1\}
$. 

Again, since $\lambda_{i+1}$ and $\lambda_{i+2}$ are electro-ambient
quasigeodesics, we further note that there exists $C > 0$ (assuming
the same $C$ for convenience) such that for all $(x,i) \in \lambda_i$,
$(x,i+1)$ lies in a (hyperbolic) $C$-neighborhood of
$\lambda_{i+1}$. Similarly for all $(x,i+2) \in \lambda_{i+2}$,
$(x,i+3)$ lies in a (hyperbolic) $C$-neighborhood of
$\lambda_{i+3}$. The same holds if we go `down' from $\lambda_{i+1}$
to $\lambda_i$ or from $\lambda_{i+3}$ to $\lambda_{i+2}$.
The {\bf
  $B_\lambda$-elementary admissible paths} in $\widetilde{B}$ consist
of the following:\\
1) Horizontal  subsegments of $\lambda_j$,  $j = \{ i,\cdots i+ 3 \}$.\\
2) Vertical segments of length $1$ joining $x  \times \{ i+1 \}$ to
$\phi (x)  \times \{ i+2 \}$, for $x \in \lambda_{i+1}$.\\
3)  Vertical segments of length $1$ joining $x  \times \{ j  \}$ to
$x  \times \{ j+1 \}$, for $j = i, i+2$.\\
4)  Horizontal geodesic segments lying in a {\em hyperbolic}
 $C$-neighborhood of
  $\lambda_j$, $j = i,\cdots i+3$.\\
5)  Horizontal hyperbolic segments of {\em electric length $\leq C$}
  and {\em hyperbolic length $\leq K(B)$} joining points of the form
  $( \phi (x), i+2)$ to a point on $\lambda_{i+2}$ for $(x, i+1) \in
  \lambda_{i+1}$. \\
6)  Horizontal hyperbolic segments of {\em electric length $\leq C$}
  and {\em hyperbolic length $\leq K(B)$} joining points of the form
  $( \phi^{-1} (x), i+1)$ to a point on $\lambda_{i+1}$ for $(x, i+2) \in
  \lambda_{i+2}$. \\
7)  Hyperbolic geodesic segments lying entirely within some lift of
  a Margulis tube $N_\epsilon ( \widetilde{\sigma} ) \times [1,2]$
  joining points $x, y \in \lambda_{i+1} \cup \lambda_{i+2}$.

\smallskip

{\bf Definition:} A   $B_\lambda$-admissible path is a union of
$B_\lambda$-elementary admissible paths. 

The following lemma follows from the above definition and Lemma
\ref{bddheight}.

\begin{lemma}
There exists a function $g: \Bbb{Z} \rightarrow \Bbb{N}$ such that
for any block $B_i$, and $x$ lying on a $B_\lambda$-admissible
path in $\widetilde{B_i}$, there exist $y
\in \lambda_{j}$ and $z \in \lambda_k$ where $\lambda_j \subset B_\lambda$ and
$\lambda_k \subset B_\lambda$ lie on the two boundary horizontal sheets, satisfying:

\begin{center}

$d(x, y) \leq g(i)$ \\
$d(x, z) \leq g(i)$ \\

\end{center}
\label{bddheight-adm}
\end{lemma}

Let $h(i) = \Sigma_{j = 0 \cdots i} g(j)$ be the sum of the values of
$g(j)$ as $j$ ranges from $0$ to $i$ (with the assumption that
increments are by $+1$ for $i \geq 0$ and by $-1$ for $i \leq
0$). Then we have from Lemma \ref{bddheight-adm} above,

\begin{cor}
There exists a function $h: \Bbb{Z} \rightarrow \Bbb{N}$ such that
for any block $B_i$, and $x$ lying on a $B_\lambda$-admissible
path in $\widetilde{B_i}$, there exist $y
\in \lambda_0 = \lambda$ such that:

\begin{center}

$d(x, y) \leq h(i)$ 

\end{center}
\label{bddht-final}
\end{cor}

{\bf Important Note:} In the above Lemma \ref{bddheight-adm} and
Corollary \ref{bddht-final}, it is important to note
that the distance $d$ is {\bf hyperbolic}, not electric. 
This is because the number $K(B_i)$ occurring in elementary paths of
type $5$ and $6$ is a hyperbolic length depending only on $i$ (in $B_i$).

\smallskip

Next suppose that $\lambda$ lies outside $B_N(p)$, the $N$-ball about
a fixed reference point $p$ on the boundary horizontal surface
$\widetilde{S} \times \{ 0 \} \subset \widetilde{B_0}$. Then by
Corollary \ref{bddht-final}, any $x$ lying on a $B_\lambda$-admissible
path in $\widetilde{B_i}$ satisfies

\begin{center}

$d(x, p) \geq N - h(i)$ 

\end{center}

Also, since the electric, and hence hyperbolic `thickness' (the
shortest distance between its boundary horizontal sheets) is $\geq 1$,
we get,

\begin{center}

$d(x, p) \geq |i|$ 

\end{center}

Assume for convenience that $i \geq 0$ (a similar argument works,
reversing signs for $i < 0$). Then,

\begin{center}

$d(x, p) \geq min \{ i,  N - h(i) \}$ 

\end{center}

Let $h_1 (i) = h(i) + i$.  Then $h_1$ is a
monotonically increasing function on the integers.
If
 $h_1^{-1} (N)$ denote the largest positive
integer $n$ such that $h(n) \leq m$, then clearly, .
 $h_1^{-1} (N) \rightarrow \infty$
 as $N \rightarrow \infty$. We have thus shown:

\begin{lemma}
There exists a function $M(N) : \Bbb{N} \rightarrow \Bbb{N}$
such that  $M (N) \rightarrow \infty$
 as $N \rightarrow \infty$ for which the following holds:\\
For any geodesic $\lambda \subset \widetilde{S} \times \{ 0 \} \subset
\widetilde{B_0}$, a fixed reference point 
 $p \in \widetilde{S} \times \{ 0 \} \subset
\widetilde{B_0}$ and any $x$ on a $B_\lambda$-admissible path, 

\begin{center}

$d(\lambda , p) \geq N \Rightarrow d(x,p) \geq M(N)$.

\end{center}

\label{far-nopunct}
\end{lemma}

{\bf Aside for Punctured Surfaces}

\smallskip

We mention parenthetically the versions  of Corollary
\ref{bddht-final} and Lemma \ref{far-nopunct} here
that will be useful for punctured surfaces in the next section.

\begin{cor}
There exists a function $h: \Bbb{Z} \rightarrow \Bbb{N}$ such that
for any block $B_i$, and $x$ lying on $\lambda_{ib}$, there exist $y
\in \lambda_{0b} = \lambda_b$ such that:

\begin{center}

$d(x, y) \leq h(i)$ 

\end{center}
\label{bddht-final-punct}
\end{cor}

\begin{lemma}
There exists a function $M(N) : \Bbb{N} \rightarrow \Bbb{N}$
such that  $M (N) \rightarrow \infty$
 as $N \rightarrow \infty$ for which the following holds:\\
For any horo-ambient quasigeodesic
 $\lambda \subset \widetilde{S} \times \{ 0 \} \subset
\widetilde{B_0}$, a fixed reference point 
 $p \in \widetilde{S} \times \{ 0 \} \subset
\widetilde{B_0}$ and any $x$ on some $\lambda_{ib}$, 

\begin{center}

$d(\lambda_b , p) \geq N \Rightarrow d(x,p) \geq M(N)$.

\end{center}

\label{far-punct}
\end{lemma}

In Lemma \ref{far-punct} we have used $\lambda_{b}$ in place of
$\lambda$ as $\lambda_b$ is constructed from $\lambda^h$ by changing
it along horocycles. However, another version of the above Lemma will
also sometimes be useful. If we start with a $\lambda^h$ that lies
outside large balls about $p$, we can ensure that $\lambda_b$ also
lies outside large balls about $p$, for $\lambda$ may approach $p$
only along cusps. Hence we may replace the hypothesis that $\lambda_b$
lie outside $B_N(p)$ by the  hypothesis that $\lambda^h$
lie outside $B_N(p)$:

\begin{lemma}
There exists a function $M(N) : \Bbb{N} \rightarrow \Bbb{N}$
such that  $M (N) \rightarrow \infty$
 as $N \rightarrow \infty$ for which the following holds:\\
For any hyperbolic geodesic
 $\lambda^h \subset \widetilde{S^h} \times \{ 0 \} \subset
\widetilde{B_0}$, a fixed reference point 
 $p \in \widetilde{S^h} \times \{ 0 \} \subset
\widetilde{B_0}$ and any $x$ on some $\lambda_{ib}$, 

\begin{center}

$d(\lambda^h , p) \geq N \Rightarrow d(x,p) \geq M(N)$.

\end{center}

\label{far-punct2}
\end{lemma}

\subsection{Joining the Dots}

Recall that {\bf admissible paths} in a model manifold of bounded
geometry consist of:\\
1) Horizontal segments along some $\widetilde{S} \times \{ i \}$
  for $ i = \{ 0,1,2,3 \}$ (thin blocks) or $i = \{ 0, 1 \}$ (thick
  blocks).\\
2) Vertical segments $x \times [0,1]$ or $x \times [2,3]$ for thin
  blocks.\\
3) Vertical segments of length $1$ joining $x  \times \{ 1 \}$ to
$\phi (x)  \times \{ 2 \}$ for thin blocks.\\
4)Vertical segments of length $1$ joining $x  \times \{ 0 \}$ to
$\phi (x)  \times \{ 1 \}$ for thick blocks.\\

\smallskip

Our strategy in this subsection is: \\
$\bullet 1$ Start with an electric geodesic $\beta_e$ in
$\widetilde{M_{el}}$ joining the end-points of $\lambda$.
\\
$\bullet 2$ Replace it by an {\em admissible quasigeodesic}, i.e. an
admissible path that is a quasigeodesic. 
\\
$\bullet 3$ Project the intersection of the
 admissible quasigeodesic with the horizontal sheets onto
 $B_\lambda$. \\
$\bullet 4$ The result of step 3 above is disconnected. {\em Join the
   dots} using $B_\lambda$-admissible paths. \\

The end product is an electric quasigeodesic built up of $B_\lambda$
admissible paths.

Now for the first two steps: Since $\widetilde{B}$ (for a thick block
$B$)  has thickness $1$, any path lying in a thick block can be
pertubed to an admissible path lying in $\widetilde{B}$, changing the
length by at most a bounded multiplicative factor. For $B$ thin, we
decompose paths into horizontal paths lying in some $\widetilde{S}
\times \{ j \}$, for $j = 0, \cdots 3$ and 
vertical paths of types (2) or (3) above. All this can be done as for
thick blocks, changing lengths by a bounded multiplicative
constant. The result is therefore an electric quasigeodesic. Without
loss of generality, we can assume that the electric quasigeodesic is
one without back-tracking (as this can be done without increasing the
length of the geodesic - see \cite{farb-relhyp} or \cite{klarreich}
for instance). Abusing notation slightly, assume therefore that 
$\beta_e$ is an admissible electric quasigeodesic without backtracking
joining the end-points of $\lambda$.

Now act on $\beta_e$ by $\Pi_\lambda$. From Theorem \ref{retract}, we
conclude, by restricting $\Pi_\lambda$ to the horizontal sheets of
$\widetilde{M_{el}}$ that the image $\Pi_\lambda ( \beta_e )$  is a
`dotted electric quasigeodesic' lying entirely on $B_\lambda$. This
completes step 3. Note that since $\beta_e$ consists of admissible
segments, we can arrange so that two nearest points on $\beta_e$ which are not
connected to each 
other are at a distance of one apart, i.e. they form the end-points of
a vertical segment of type (2), (3) or (4). Let $\Pi_\lambda ( \beta_e
) \cap B_\lambda = \beta_d$, be the dotted quasigedoesic lying on
$B_\lambda$. We want to join the dots in $\beta_d$ converting it into
a {\bf connected} electric quasigeodesic built up of {\bf
  $B_\lambda$-admissible paths}.  

For vertical segments of type (4) joining $p, q$ (say), $\Pi_\lambda
(p), \Pi_\lambda (q)$ are a bounded hyperbolic distance apart. Hence, by
the proof of Lemma \ref{retract-thick}, we can join 
 $\Pi_\lambda
(p), \Pi_\lambda (q)$ by a $B_\lambda$-admissible path of length
bounded by some $C_0$ (independent of $B$, $\lambda$).

For vertical segments of type (2) joining $p, q$, we note that $\Pi_\lambda
(p), \Pi_\lambda (q)$ are a bounded hyperbolic distance apart. Hence, by
the proof of Lemma \ref{retract-thin}, we can join 
 $\Pi_\lambda
(p), \Pi_\lambda (q)$ by a $B_\lambda$-admissible path of length
bounded by some $C_1$ (independent of $B$, $\lambda$).

This leaves us to deal with case (3). Such a  segment consists of
a segment lying within a lift of a
Margulis tube and a horizontal segment of length $1$ lying
outside. Decompose the bit within a Margulis tube into a horizontal
segment lying on some horizontal surface and (possibly) a vertical
segment of hyperbolic length $1$. The image of the horizontal
part of the  path  is again uniformly bounded in the electric metric.
Further, by {\bf co-boundedness}, we can ensure that the {\em hyperbolic}
length of the image away from the lift of 
at most one Margulis tube lying in the zero neighborhood of
$B_\lambda$
is bounded uniformly by some $C_2$. The same can be ensured of
vertical paths of hyperbolic length one lying inside lifts of 
Margulis tubes. These two pieces (images under $\Pi_\lambda$
of horizontal paths inside
lifts of Margulis tubes and vertical segments of length one inside
lifts of Margulis tubes) can be replaced by $B_\lambda$-admissible
paths of uniformly bounded electric length (since at most one lift of
a Margulis tube lying in a zero neighborhood of $B_\lambda$ is
in the image for length $\geq C_2$.) Finally, the segment lying
outside, being horizontal, its image is connected and of bounded
length by Lemma \ref{easyprojnlemma}.

After joining the dots, we can assume further that the quasigeodesic
thus obtained does not backtrack (cf \cite{farb-relhyp} and
\cite{klarreich}).

Putting all this together, we conclude:

\begin{lemma}
There exists a function $M(N) : \Bbb{N} \rightarrow \Bbb{N}$
such that  $M (N) \rightarrow \infty$
 as $N \rightarrow \infty$ for which the following holds:\\
For any geodesic $\lambda \subset \widetilde{S} \times \{ 0 \} \subset
\widetilde{B_0}$, and  a fixed reference point 
 $p \in \widetilde{S} \times \{ 0 \} \subset
\widetilde{B_0}$,
there exists a connected electric quasigeodesic $\beta_{adm}$  without 
backtracking, such that \\
$\bullet$ $\beta_{adm}$ is built up of $B_\lambda$-admissible
paths. \\
$\bullet$ $\beta_{adm}$ joins
the end-points of $\lambda$. \\
$\bullet$ 
$d(\lambda , p) \geq N \Rightarrow d(\beta_{adm},p) \geq M(N)$. \\
\label{adm-qgeod-props}
\end{lemma}

{\bf Proof:} The first two criteria follow from the discussion
preceding this lemma. The last follows from Lemma \ref{far-nopunct}
since the discussion above gives a quasigeodesic built up out of
admissible paths.
$\Box$

\subsection{Proof of Theorem}

{ \bf Electric Geometry Revisited}

We note the following properties of the pair $(X, \mathcal{H})$
where $X = \widetilde{M}$ and $\mathcal{H}$ consists of the lifts of
Margulis tubes in $M$ to
the universal cover. Each lift of a Margulis tube shall henceforth be
termed an {\em extended Margulis tube}. There exist $C, D, \Delta$ such that\\
1) Each extended Margulis tube is $C$-quasiconvex.\\
2)  Any two extended Margulis tubes are $D$-separated.\\
3)  The collection $\mathcal{H}$ is $C$-cobounded, i.e. the nearest
  point projection of any member of $\mathcal{H}$ onto any other has
  diameter bounded by $C$.\\
4) $\widetilde{M_{el}} = X_{el}$ is $\Delta$-hyperbolic, (where
  $\widetilde{M_{el}} = X_{el}$ is the electric metric on
  $\widetilde{M} = X$
  obtained by electrocuting all extended Margulis tubes, i.e. all
  members of $\mathcal{H}$).\\
5) $({X_{el}}, \mathcal{H} )$ has the {\em Bounded Penetration
  Property}.\\
6)  An electro-ambient quasigeodesic is a hyperbolic quasigeodesic.\\

\smallskip

The first property follows from the fact that each
$\epsilon$-neighborhood of a closed geodesic in a hyperbolic manifold
is convex for sufficiently small $\epsilon$. 

The second follows from choosing $\epsilon$ sufficiently small.

The third follows from the uniform separation of (the convex) extended
Margulis sets. 

The fourth and fifth follow from Lemmas \ref{farb1A} and \ref{farb2A}.

Property (6) follows from Lemma \ref{ea2}.

{\bf Note:} So far we have not used the hypothesis that
$\widetilde{M}$ and hence, (from Property (4) above, or by Lemma
\ref{farb1A})
 $\widetilde{M_{el}}$ are hyperbolic metric
spaces. It is at this stage that we shall do so and assemble the proof
of the main Theorem.

\begin{theorem}
Let $M$ be a 3 manifold homeomorphic to $S \times J$ (for $J = [0,
  \infty ) $ or $( - \infty , \infty )$). Further suppose that $M$ has
  {\em i-bounded geometry}, where $S_0 \subset B_0$ is the lower
  horizontal surface of the building block $B_0$. Then the inclusion
  $i : \widetilde{S} \rightarrow \widetilde{M}$ extends continuously
  to a map 
  $\hat{i} : \widehat{S} \rightarrow \widehat{M}$. Hence the limit set
  of $\widetilde{S}$ is locally connected.
\label{crucial}
\end{theorem}

{\bf Proof:}
Suppose $\lambda \subset \widetilde{S}$ lies outside
a large $N$-ball about $p$. By Lemma 
\ref{adm-qgeod-props} we obtain an electric quasigeodesic without backtracking
$\beta_{adm}$
built up of
$B_\lambda$-admissible  paths
 lying outside an $M(N)$-ball about $p$ (where $M(N) \rightarrow
 \infty $ as $N \rightarrow \infty $). 

Suppose that $\beta_{adm}$ is
 a $(K, \epsilon)$ quasigeodesic. Note that $K, \epsilon$ depend on
 `the Lipschitz constant' of $\Pi_\lambda$ and hence only on
 $\widetilde{S}$ and $\widetilde{M}$.

From Property (6) (or Lemma \ref{ea-genl}) we find that
if $\beta^{h}$ denote the hyperbolic geodesic in
$\widetilde{M}$ joining the end-points of $\lambda$, then $\beta^h$
lies in a (uniform) $C^{\prime}$ neighborhood of $\beta_{adm}$. 

Let $M_1(N) = M(N) - C^{\prime}$.
Then $M_1(N) \rightarrow
\infty$ as $N \rightarrow \infty$. Further, the hyperbolic geodesic 
 $\beta^h$ 
 lies outside an $M_1(N)$-ball around $p$. Hence, by Lemma
 \ref{contlemma}, 
the inclusion
  $i : \widetilde{S} \rightarrow \widetilde{M}$ extends continuously
  to a map 
  $\hat{i} : \widehat{S} \rightarrow \widehat{M}$. 

Since the
  continuous image of a compact locally connected set is locally
  connected (see \cite{hock-young} )
and the (intrinsic) boundary of $\widetilde{S}$ is a circle, we
  conclude that the limit set
  of $\widetilde{S}$ is locally connected.

This proves the theorem. $\Box$

\section{Cannon-Thurston Maps for Surfaces With Punctures}

\subsection{Modification of Construction for Punctured Surfaces}

\smallskip

We summarise the modifications to be made to the construction in the
previous sections, so as to make the results applicable for punctured
surfaces:

\begin{enumerate}

\item $\lambda$ is now a horo-ambient quasigeodesic built out of a
  hyperbolic geodesic $\lambda^h$ \\
\item $\Pi_\lambda$ and $B_\lambda$ are constructed as before \\
\item Let $\beta_a$ be an ambient admissible quasigeodesic, i.e. an
  ambient quasigeodesic built up of elementary admissible paths. Then 
$\Pi_\lambda ( \beta_a \cap \widetilde{M_{BH}} ) \subset
  B_\lambda$. \\
\item Joining the dots on this projected image of $\beta_a$, 
we obtain finally via Lemma \ref{far-punct} \\
$\bullet$ Suppose $\lambda^h$ lies outside large balls about a fixed
reference point $p$. There exists an ambient admissible electric quasigeodesic
$\beta_{amb}$ in $\widetilde{M_{el}}$ such that any horizontal piece
of $\beta_{amb} \cap B^b_\lambda$ also lies outside large
balls. Further, any piece of $\beta_{amb}$ lying wholly inside the
lift of a
Margulis tube also lies outside large balls. (Note that $B^b_\lambda =
\bigcup_i \lambda_{ib}$). 

\end{enumerate}

Next recall 
that $M^h$ is a hyperbolic manifold with  cusps, and that excising
these cusps we get $M$ which is naturally homeomorphic to $S \times J$
where 
$J = \Bbb{R}$ 
or $[0 , \infty )$. Further we may assume that $M$ is given the
  structure of a model of {\bf i-bounded geometry}. (Here we are
  abusing notation slightly as $M$, strictly speaking, is only
  quasi-isometric to a model of i-bounded geometry.) Let $M^h_{el}$, $M_{el}$
  denote $M^h$, $M$ with Margulis tubes electrocuted. 
Let ${\mathcal{H}}_0$ denote the collection of horoballs that
corespond to the lifts of cusps in $\widetilde {M^h_{el}}$. Thus, 
 $\widetilde{M_{el}} =  \widetilde {M^h_{el}} - \{ H : H \in
{\mathcal{H}}_0 \}$. 
$ \widetilde {M^h_{el}}$ is hyperbolic by Lemma \ref{farb1A}.
\smallskip

Now let $\beta^h_{el}$ be the electric geodesic in the hyperbolic
metric space $\widetilde{M^h_{el}}$ joining the end-points $a, b$ of
$\lambda^h$. Let $H(\beta^h_{el} )$ denote the union of $\beta^h_{el}$
and the collection of horoballs in ${\mathcal{H}}_0$ that
$\beta^h_{el}$ meets. Then by 
Theorem \ref{ctm} (using the fact stated there that the theroem goes
through for separated mutually cobounded uniformly quasiconvex sets),
we have \\ 

\noindent $\bullet 1$ $H(\beta^h_{el} )$ is quasiconvex in
$\widetilde{M^h_{el}}$. \\
\noindent $\bullet 2$ $\beta_{amb}$ lies in a bounded electric neighborhood of 
$H(\beta^h_{el} )$ \\

\subsection{Electrically close implies hyperbolically close}

In what follows we want to construct out of $\beta_{amb}$ a
hyperbolic quasigeodesic $\gamma$ in $\widetilde {M^h_{el}}$ such that
entry and exit points of $\gamma$ with respect to $H \in
{\mathcal{H}}_0$ lie outside large balls $B_N (p) \subset
\widetilde{M}$ (here the metric is the {\em hyperbolic} metric). The
strategy is as follows: \\

For any $H_i \in {\mathcal{H}}_0$ look at the part $\beta_i$ of
$\beta_{amb}$ that lies close to $H_i$. We claim that if this piece is
long, then after pruning it a bit at the ends if necessary, the pruned
subsegment of $\beta_i$ lies {\bf hyperbolically close} to $H_i$. We
make this precise below. 

By Theorem \ref{ctm}, and as in \cite{ctm-locconn},
 there exists $\Delta  > 0$ such that
$\beta_{amb}$ lies in an (electric) $\Delta$ neighborhood of $H(
\beta^h_{el})$. Let $H_1, \cdots H_k$ denote the horoballs in 
$H(
\beta^h_{el})$. Let $\beta_i$ be the maximal subsegment of
$\beta_{amb}$ joining points of $N_\Delta^{el} (H_i) \cap
\beta_{amb}$. Then there exists $D = D( \Delta )$ such that $ \beta_i
\subset N_D^{el} ( H_i )$. Let $a_i, b_i$ be the end-points of
$\beta_i$ and $P_i$ denote nearest point projection onto $H_i$. 
$[x,y]_e$ will denote the {\em electric} geodesic joining $x, y$.
$[P_i(x), P_i(y)]$ will denote the {\em hyperbolic} geodesic joining
 $P_i(x), P_i(y)$ within the horoball $H_i$. 

Fixing $K \geq 0$ ($K = 4D$ will suffice for our purposes) let $c_i,
d_i \in \beta_i$ be such that $\overline{c_i d_i}$, the subsegment of
$\beta_i$ joining $c_i,d_i$  has length less
than $K$. Suppose further that $a_i, c_i, d_i, b_i$ occur in that
order along the segment joining $a_i, b_i$. 
Then \\ $[P_i (c_i), c_i]_e \cup \overline{c_i d_i} \cup [d_i,
  P_i(d_i)]_e \cup [P_i(d_i),P_i(c_i)] = \sigma$ is a loop of electric
length
less than $C = C(K,D)$ ($ = C(D)$ if $K = 4C$). This follows from the
following observations:\\
1) $[P(c_i), c_i]_e$,  $[P(d_i), d_i]_e$ have length less than or
equal to $D$ \\
2) $\overline{c_i d_i}$  has length less
than $K$ \\
3) $[P_i(c_i),P_i(d_i)]$ has length bounded in terms of $K$ by Lemma
\ref{easyprojnlemma}. \\

Since $\sigma$ has electric length less than $C $, we could conclude
that $\sigma$ has bounded hyperbolic length by Lemma \ref{loop-bdd} if
in addition we knew that $\sigma$ does not backtrack. (In particular we would
be able to show that $\sigma$ has bounded penetration.) However, we
only know that each of the four components of $\sigma$ individually
does not backtrack. In fact,  
$[ c_i, P_i (c_i)]_e  \cup [P_i(c_i),P_i(d_i)] \cup [P_i(d_i),d_i]_e
$ is a path without backtracking. Therefore, backtracking, if it exists, is  a 
consequence of overlap of initial segments of $[c_i.P_i(c_i)]_e$ and
$\overline{c_i d_i}$ (or, $[d_i.P_i(d_i)]_e$ and
$\overline{d_i, c_i}$). Clearly, such overlaps can have length at most
$D$. Therefore, any such segment $\overline{c_i d_i}$  with $d_e(a_i,
c_i) \geq D$, $d_e(b_i,
d_i) \geq D$ must have bounded penetration property (since the paths
$[ a_i, P_i (a_i)]_e  \cup [P_i(a_i),P_i(b_i)] \cup [P_i(b_i),b_i]_e
$ and $\overline{a_i b_i}$ can have overlaps of length at most $D$ at
the beginning and end), i.e there exists
$D_0 = D_0(D,K) \geq 0$ such that $\overline{c_i d_i} \cap T$ has {\em
  hyperbolic 
  length} less than $D_0$ (where $T$ is any
lift of a Margulis tube). 

Now, choose $x \in \overline{a_i b_i}$ such that $d_e(x,a_i) \geq 2D$
and  $d_e(x,b_i) \geq 2D$. Choose $c_i, d_i$ such that  $d_e(x,c_i) =
2D$, 
 $d_e(x,d_i) = 2D$ and $a_i, c_i, d_i, b_i$ lie in that order along
the path from $a_i$ to $b_i$. Then using the loop
$[P_i (c_i), c_i]_e \cup \overline{c_i,x} \cup [x,
  P_i(x)]_e \cup [P_i(x),P_i(c_i)] = \sigma$ and the argument above,
we conclude that $[x,P(x)]_e$ satisfies the {\em bounded penetration}
property and hence $[x,P(x)]$ has bounded hyperbolic length. This is summarised
in the following Lemma.

\begin{lemma}
There exists $D_0 \geq 0$ such that the following holds. Let $\beta_i
= \overline{a_i b_i}$ be as above and $x \in \overline{a_i b_i}$ with 
 $d_e(x,a_i) \geq
2D$, 
 $d_e(x,b_i) \geq 2D$. Then $d(x,H_i) \leq D_0$. (Note that $d(x,H_i)$
denotes hyperbolic distance.)
\label{hypclose}
\end{lemma}

Thus the subpath of $\beta_i$ obtained by pruning pieces of (electric) length
$2D$ from the beginning and the end lies in a bounded {\em hyperbolic}
neighborhood of $H_i$ (and not just in a bounded electric neighborhood
of $H_i$). 

\subsection{Constructing an electric quasigeodesic}

The argument in this subsection is a slight modification of the
argument in \cite{brahma-pared} for punctured surfaces of bounded
geometry. The slight modification is due to  pruning 
electric quasigeodesics that follow a horoball for a considerable length.

Choose from the the collection of $H_i \in \mathcal{H} (\beta )$ the
subcollection for which $\beta_i$ has diameter greater than $4D$. We
denote this subcollection as ${\mathcal{H}}_l (\beta )$ ($l$ stand for
`large'). Let $H_{l1}, \cdots H_{lk}$ be the horoballs in this
collection. 

For the relevant subpaths $\beta_{l1}, \cdots \beta_{lk}$ of $\beta$
we construct $\gamma_{l1}, \cdots \gamma_{lk}$ as follows.

Let $\alpha_{li} = \overline{c_{li} d_{li}} \subset \beta_{li} =
\overline{a_{li} b_{li}}$ denote the 
subpath at distance less than or equal to $D$ from $H_{li}$. By Lemma
\ref{hypclose} we have $d_e(c_{li}, a_{li}) \leq 2D$ and 
 $d_e(d_{li}, b_{li}) \leq 2D$. Let \\
$\gamma_{li} = 
[c_{li},P_{li} (c_{li})]_e \cup  [P_{li}(c_{li}),P_{li}(d_{li})] \cup
[P_{li}(d_{li}), d_{li}]_e $ \\
Let $\gamma = ( \beta - \bigcup_i \beta_{li} ) \cup \bigcup_i
\gamma_{li}$.

Each $\beta_{li}$ starts and ends (electrically) close to the entry
and exit points 
$u_{li},v_{li}$ of $\beta^h_{el}$ with respect to the horoball
$H_{li}$. 

Since $c_{li}, d_{li}$ are close (bounded by $2D$) to
$a_{li}, b_{li}$ respectively, then from Lemma
\ref{coboundedHor&T}, Lemma \ref{loop-bdd} and Lemma \ref{easyprojnlemma}
we find that there exists $D_1 \geq 0$ such that \\
$d(P_{li}(c_{li}),P_{li}(a_{li})) \leq D_1$ \\
$d(P_{li}(d_{li}),P_{li}(b_{li})) \leq D_1$ \\

Note that $d$ here is the {\em hyperbolic} distance. Hence the
hyperbolic geodesic $[c_{li},d_{li}]$ lies close to  $[a_{li},b_{li}]$
and hence to the hyperbolic geodesic  $[u_{li},v_{li}]$ (by fellow
traveller property).

Thus we conclude \\
$\bullet$ $\gamma$ lies in a bounded neighbourhood of the electric
geodesic $\beta^h_{el}$. 

\smallskip

{\bf Note:} The remaining $\beta_i$'s being less than $4D$ in
length are therefore uniformly bounded. Hence their projections onto
the corresponding $H_i$'s are also uniformly bounded in diameter. The
length of $\beta^h_{el} \cap H_i$ for these $H_i$'s is also therefore
uniformly bounded. (Else the projection onto $\beta^h_{el} \cup
\mathcal{H} ({\beta^h_{el}})$ would have to have jumps and hence not
be `large-scale continuous'.)

\smallskip

Since $\gamma$ is obtained from $\beta_{amb}$, $\gamma$ tracks
$\beta^h_{el}$ off horoballs. Further, since entry and exit points of
$\gamma$ and $\beta^h_{el}$ with respect to horoballs are a bounded
distance apart, they are fellow travellers within horoballs. 
From this it follows easily that $\gamma$ is an electric
quasigeodesic. 

$\gamma$ therefore has two properties: \\
1) $\gamma$ lies close to $\beta^h_{el}$ and is an electric
  quasigeodesic. \\
2) All points of $\gamma \cap \widetilde{M}$ lie outside a large ball
  about the fixed reference point $p$ if $\lambda^h$ does. 

\smallskip

The first property follows from the above discussion and the last is
just a restatement of property (4) of Section 7.1 (the first
subsection of the present section), coupled with the fact that entry
points of $\gamma$ into horoballs $H_i \in \mathcal{H}$ lie {\em
  hyperbolically} close to $\beta_{amb}$.

This gives rise to the following property of $\gamma$. Recall that
building blocks are built from the truncated surface $S$, and that we
fix a `starting block' $B_0$. We identify $S \times \{ 0 \}$ with the
truncated surface obtained from $S^h$ the hyperbolic reference
surface. 

\begin{prop}
There exists a function $M(N) : \Bbb{N} \rightarrow \Bbb{N}$
such that  $M (N) \rightarrow \infty$
 as $N \rightarrow \infty$ for which the following holds:\\
For any geodesic $\lambda^h \subset \widetilde{S^h}$, 
and  a fixed reference point 
 $p \in \widetilde{S} \times \{ 0 \} \subset
\widetilde{B_0}$,
there exists a connected electric quasigeodesic $\gamma$  without 
backtracking, such that \\
$\bullet$ If $\lambda^h$ lies outside $B_N(p)$, then every point $x$ of
$\gamma - \{ H : H \in \mathcal{H} \}$ lies at a {\em hyperbolic}
distance of at
least $M(N)$ from $p$.
 \\
\label{farOhoroballs}
\end{prop}

The above Proposition is a punctured surface version of 
Lemma \ref{adm-qgeod-props}.

Now, recall a Lemma from \cite{brahma-pared} (which has been proven
there as a part of Theorem 5.9).

\begin{lemma}
There exists a function $M_1(N)$ such that
$M_1(N) \rightarrow \infty$ as $N \rightarrow \infty$ satisfying the
following: \\
Given a uniformly separated collection of horoballs $\mathcal{H}$ and a point
$p$ lying outside them, let $\gamma$ be a path without backtracking,
such that $\gamma - \{ H : H \in \mathcal{H} \}$ lies
outside $B_N(p)$. Further suppose that $\gamma \cap H$ is a
(hyperbolic) geodesic, whenever  $\gamma \cap H$ is non-empty.
 Then $\gamma$ lies outside an $M_1(N)$ ball about $p$.
\label{horo-far}
\end{lemma}

Combining Proposition \ref{farOhoroballs} (for pieces of $\gamma$
outside horoballs)  and Lemma \ref{horo-far}
above  (for the geodesic
segments within horoballs) we conclude: \\

\begin{prop}
There exists a function $M(N) : \Bbb{N} \rightarrow \Bbb{N}$
such that  $M (N) \rightarrow \infty$
 as $N \rightarrow \infty$ for which the following holds:\\
For any geodesic $\lambda^h \subset \widetilde{S^h}$, 
and  a fixed reference point 
 $p \in \widetilde{S} \times \{ 0 \} \subset
\widetilde{B_0}$,
there exists a connected electric quasigeodesic $\gamma$  without 
backtracking, such that \\
$\bullet$ If $\lambda^h$ lies outside $B_N(p)$, then every point $x$ of
$\gamma$ lies at a {\em hyperbolic} distance of at
least $M(N)$ from $p$.
\label{farforpunct}
\end{prop}

\subsection{From electric quasigeodesics to hyperbolic quasigeodesics}

We have thus constructed an electric quasigeodesic $\gamma$ without
backtracking joining the end-points of $\lambda^h$ every point of
which lies outside a (hyperbolic) large ball about $p$. The last step
is to promote $\gamma$ to a {\em hyperbolic quasigeodesic}.

Since $\gamma$ is built up of admissible paths within Margulis tubes,
we might as well assume that $\gamma$ is an electro-ambient
quasigeodesic without backtracking. 

\begin{lemma}
The undelying path of $\gamma$ is a hyperbolic quasigeodesic.
\label{gamma-qg}
\end{lemma}

{\bf Proof:} Margulis tubes satisfy the mutual co-boundedness property
by Lemma \ref{coboundedHor&T}. Hence by Lemma \ref{ea2}, $\gamma$ is a
hyperbolic quasigeodesic. $\Box$

\begin{theorem}
Let $M^h$ be a 3 manifold homeomorphic to $S^h \times J$ (for $J = [0,
  \infty ) $ or $( - \infty , \infty )$). Further suppose that $M^h$ has
  {\em i-bounded geometry}. Let $S_0 \subset B_0$ be the lower
  horizontal surface of the building block $B_0$ in the manifold
$M$  obtained by removing cusps. Then the inclusion
  $i : \widetilde{S^h} \rightarrow \widetilde{M^h}$ extends continuously
  to a map 
  $\hat{i} : \widehat{S^h} \rightarrow \widehat{M^h}$. Hence the limit set
  of $\widetilde{S}$ is locally connected.
\label{crucial-punct}
\end{theorem}

{\bf Proof:}
Suppose $\lambda^h \subset \widetilde{S^h}$ lies outside
a large $N$-ball about $p$. By Lemma 
\ref{farforpunct} and Lemma \ref{gamma-qg}, 
 we obtain a hyperbolic quasigeodesic $\gamma$
 lying outside an $M(N)$-ball about $p$ (where $M(N) \rightarrow
 \infty $ as $N \rightarrow \infty $). 

If $\beta^{h}$ denote the hyperbolic geodesic in
$\widetilde{M^h}$ joining the end-points of $\lambda^h$, then $\beta^h$
lies in a (uniform) $C^{\prime}$ neighborhood of $\gamma$ (since
hyperbolic quasigeodesics starting and ending at the same points track
each other throughout their lengths). 

Let $M_1(N) = M(N) - C^{\prime}$.
Then $M_1(N) \rightarrow
\infty$ as $N \rightarrow \infty$. Further, the hyperbolic geodesic 
 $\beta^h$ 
 lies outside an $M_1(N)$-ball around $p$. Hence, by Lemma
 \ref{contlemma}, 
the inclusion
  $i : \widetilde{S^h} \rightarrow \widetilde{M^h}$ extends continuously
  to a map 
  $\hat{i} : \widehat{S^h} \rightarrow \widehat{M^h}$. 

Since the
  continuous image of a compact locally connected set is locally
  connected (see \cite{hock-young} )
and the (intrinsic) boundary of $\widetilde{S^h}$ is a circle, we
  conclude that the limit set
  of $\widetilde{S^h}$ is locally connected.

This proves the theorem. $\Box$

The proof of the above theorem is just a modification of Theorem
\ref{crucial}, once Lemma 
\ref{farforpunct} and Lemma \ref{gamma-qg} are in place.

\bibliography{i-bdd}
\bibliographystyle{alpha}

\end{document}